\input amstex\documentstyle{amsppt}  
\pagewidth{12.5cm}\pageheight{19cm}\magnification\magstep1
\topmatter
\title A bar operator for involutions in a Coxeter group\endtitle
\author G. Lusztig\endauthor
\address{Department of Mathematics, M.I.T., Cambridge, MA 02139}\endaddress
\thanks{Supported in part by National Science Foundation grant DMS-0758262.}\endthanks
\endtopmatter   
\document
\redefine\pm{\sigma}

\define\huM{\hat{\uM}}

\define\uca{\un{\ca}}
\define\ufH{\un{\fH}}

\define\uM{\un M}

\define\ha{\hat a}

\define\si{\sim}

\define\ovs{\overset}

\define\lb{\linebreak}

\define\op{\oplus}
   
\redefine\sp{\spadesuit}
\define\part{\partial}
\define\em{\emptyset}

\define\m{\mapsto}
\define\do{\dots}

\define\sm{\smallmatrix}
\define\esm{\endsmallmatrix}
\define\sub{\subset}    

\define\T{\times}
\define\ti{\tilde}
\define\nl{\newline}
\redefine\i{^{-1}}
\define\fra{\frac}
\define\un{\underline}
\define\ov{\overline}
\define\ot{\otimes}

\define\Hom{\text{\rm Hom}}

\define\di{\diamond}

\define\a{\alpha}
\redefine\b{\beta}

\define\g{\gamma}
\redefine\d{\delta}
\define\e{\epsilon}
\define\et{\eta}

\define\p{\pi}
\define\ph{\phi}
\define\ps{\psi}
\define\r{\rho}
\define\s{\sigma}
\redefine\t{\tau}

\define\k{\kappa}
\redefine\l{\lambda}
\define\z{\zeta}
\define\x{\xi}

\define\Om{\Omega}

\define\Ph{\Phi}

\redefine\ss{\bold s}
\redefine\tt{\bold t}

\define\CC{\bold C}

\define\II{\bold I}

\define\NN{\bold N}

\define\PP{\bold P}
\define\QQ{\bold Q}

\define\SS{\bold S}

\define\ZZ{\bold Z}

\define\ca{\Cal A}

\define\ch{\Cal H}

\define\cm{\Cal M}

\define\car{\Cal R}

\define\ct{\Cal T}

\define\fH{\frak H}

\define\tb{\ti b}

\define\bul{\bullet}

\define\KL{KL}
\define\LSI{L1}
\define\HEC{L2}
\define\LV{LV}
\define\VO{V}
\define\KI{Ki}
\define\dM{\dot M}
\define\dfH{\dot{\fH}}
\define\ocT{\ovs\circ\to T}

\head Introduction and statement of results\endhead
\subhead 0.0\endsubhead
In \cite{\LV} it was shown that the vector space spanned by the involutions in a Weyl group carries
a natural Hecke algebra action and a certain bar operator. These were used in \cite{\LV} to construct a new basis
of that vector space, in the spirit of \cite{\KL}, and to give a refinement of the polynomials $P_{y,w}$ of 
\cite{\KL} in the case where $y,w$ were involutions in the Weyl group in the sense that $P_{y,w}$ was split
canonically as a sum of two polynomials with cofficients in $\NN$. 
 However, the construction of the Hecke algebra action and that
of the bar operator, although stated in elementary terms, were established in a non-elementary way. (For example,
the construction of the bar operator in \cite{\LV} was done using ideas from geometry such as Verdier duality
for $l$-adic sheaves.) In the present paper we construct the Hecke algebra action and the bar 
operator in an entirely elementary way, in the context of arbitrary Coxeter groups.

Let $W$ be a Coxeter group with set of simple reflections denoted by $S$. Let $l:W@>>>\NN$ be the standard 
length function. For $x\in W$ we set $\e_x=(-1)^{l(x)}$. Let $\le$ be the Bruhat order on $W$. Let $w\m w^*$ be 
an automorphism of $W$ with square $1$ which leaves $S$ stable, so that $l(w^*)=l(w)$ for any $w\in W$. Let 
$\II_*=\{w\in W;w^{*-1}=w\}$. (We write $w^{*-1}$ instead of $(w^*)\i$.) The elements of $\II_*$ are said 
to be {\it $*$-twisted involutions} of $W$.

Let $u$ be an indeterminate and let $\ca=\ZZ[u,u\i]$. Let $\fH$ be the free $\ca$-module with basis 
$(T_w)_{w\in W}$ with the unique $\ca$-algebra structure with unit $T_1$ such that 

(i) $T_wT_{w'}=T_{ww'}$ if $l(ww')=l(w)+l(w')$ and

(ii) $(T_s+1)(T_s-u^2)=0$ for all $s\in S$.
\nl
This is an Iwahori-Hecke algebra. (In \cite{\LV}, the notation $\fH'$ is used instead of $\fH$.)

Let $M$ be the free $\ca$-module with basis $\{a_w;w\in\II_*\}$. We have the following result which, in the 
special case where $W$ is a Weyl group or an affine Weyl group, was proved in \cite{\LV} (the general case was 
stated there without proof).

\proclaim{Theorem 0.1}There is a unique $\fH$-module structure on $M$ such that for any $s\in S$ and any 
$w\in\II_*$ we have

(i) $T_sa_w=ua_w+(u+1)a_{sw}$ if $sw=ws^*>w$;

(ii) $T_sa_w=(u^2-u-1)a_w+(u^2-u)a_{sw}$ if $sw=ws^*<w$;

(iii) $T_sa_w=a_{sws^*}$ if $sw\ne ws^*>w$;

(iv) $T_sa_w=(u^2-1)a_w+u^2a_{sws^*}$ if $sw\ne ws^*<w$.
\endproclaim
The proof is given in \S2 after some preparation in \S1.

Let $\,\bar{}\,:\fH@>>>\fH$ be the unique ring involution such that $\ov{u^nT_x}=u^{-n}T_{x\i}\i$ for any
$x\in W,n\in\ZZ$ (see \cite{\KL}). We have the following result.

\proclaim{Theorem 0.2}(a) There exists a unique $\ZZ$-linear map $\,\bar{}\,:M@>>>M$ such that 
$\ov{hm}=\bar{h}\bar{m}$ for all $h\in\fH,m\in M$ and $\ov{a_1}=a_1$. For any $m\in M$ we have $\ov{\ov{m}}=m$. 

(b) For any $w\in\II_*$ we have $\ov{a_w}=\e_wT_{w\i}\i a_{w\i}$.
\endproclaim
The proof is given in \S3. Note that (a) was conjectured in \cite{\LV} and proved there in the special case where
$W$ is a Weyl group or an affine Weyl group; (b) is new even when $W$ is a Weyl group or affine Weyl group.

\subhead 0.3\endsubhead
Let $\uca=\ZZ[v,v\i]$ where $v$ is an indeterminate. We view $\ca$ as a subring of $\uca$ by setting $u=v^2$. Let
$\uM=\uca\ot_\ca M$. We can view $M$ as an $\ca$-submodule of $\uM$. We extend $\,\bar{}\,:M@>>>M$ to a 
$\ZZ$-linear map $\,\bar{}\,:\uM@>>>\uM$ in such a way that $\ov{v^nm}=v^{-n}\ov{m}$ for $m\in M,n\in\ZZ$. For 
each $w\in\II_*$ we set $a'_w=v^{-l(w)}a_w\in\uM$. Note that $\{a'_w;w\in\II_*\}$ is an $\uca$-basis of $\uM$. 
Let $\uca_{\le0}=\ZZ[v\i]$, $\uca_{<0}=v\i\ZZ[v\i]$, $\uM_{\le0}=\sum_{w\in\II_*}\uca_{\le0}a'_w\sub\uM$,
$\uM_{<0}=\sum_{w\in\II_*}\uca_{<0}a'_w\sub\uM$.

Let $\ufH=\uca\ot_\ca\fH$. This is naturally an $\uca$-algebra containing $\fH$ as an $\ca$-subalgebra. Note that
the $\fH$-module structure on $M$ extends by $\uca$-linearity to an $\ufH$-module structure on $\uM$. We denote 
by $\,\bar{}:\uca@>>>\uca$ the ring involution such that $\ov{v^n}=v^{-n}$ for $n\in\ZZ$. We denote by 
$\,\bar{}:\ufH@>>>\ufH$ the ring involution such that $\ov{v^nT_x}=v^{-n}T_{x\i}\i$ for $n\in\ZZ,x\in W$. We have
the following result which in the special case where $W$ is a Weyl group or an affine Weyl group the theorem is 
proved in \cite{\LV, 0.3}.

\proclaim{Theorem 0.4}(a) For any $w\in\II_*$ there is a unique element 
$$A_w=v^{-l(w)}\sum_{y\in\II_*;y\le w}P^\pm_{y,w}a_y\in\uM$$
($P^\pm_{y,w}\in\ZZ[u]$) such that $\ov{A_w}=A_w$, $P^\pm_{w,w}=1$ and for any $y\in\II_*$, $y<w$, we have 
$\deg P^\pm_{y,w}\le(l(w)-l(y)-1)/2$. 

(b) The elements $A_w$ ($w\in\II_*$) form an $\uca$-basis of $\uM$.
\endproclaim
The proof is given in \S4. 

\subhead 0.5\endsubhead
As an application of our study of the bar operator we give (in 4.7) an explicit description of the M\"obius 
function of the partially ordered set $(\II_*,\le)$; we show that it has values in $\{1,-1\}$. This description 
of the M\"obius function is used to show that the constant term of $P_{y,w}^\pm$ is $1$, see 4.10.
In \S5 we study the "$K$-spherical" submodule $\uM^K$ of $\uM$ (where $K$ is a subset of $S$ which generates a
finite subgroup $W_K$ of $S$). In 5.6(f) we show that $\uM^K$ contains any element $A_w$ where $w\in\II_*$ has
maximal length in $W_KwW_{K^*}$. This result is used in \S6 to describe the action of $u\i(T_s+1)$ (with $s\in S$)
in the basis $(A_w)$ by supplying an elementary substitute for a geometric argument in \cite{\LV}, see Theorem 6.3
which was proved earlier in \cite{\LV} for the case where $W$ is a Weyl group. In 7.7 we give an inversion 
formula for the polynomials $P^\pm_{y,w}$ (for finite $W$) which involves the M\"obius function above and the 
polynomials analogous to $P^\pm_{y,w}$ with $*$ replaced by its composition with the opposition automorphism of 
$W$. In \S8 we formulate a conjecture (see 8.4) relating $P^\pm_{y,w}$ for certain twisted involutions $y,w$ in an
affine Weyl group to the $q$-analogues of weight multiplicities in \cite{\LSI}. In \S9 we show that 
for $y\le w$ in $\II_*$, $P^\pm_{y,w}$ is equal to the polynomial $P_{y,w}$ of \cite{\KL} plus an element in 
$2\ZZ[u]$. This follows from \cite{\LV} in the case where $W$ is a Weyl group. 

\subhead 0.6\endsubhead {\it Notation.} 
If $\Pi$ is a property we set $\d_\Pi=1$ if $\Pi$ is true and $\d_\Pi=0$ if $\Pi$ is false. We write $\d_{x,y}$ 
instead of $\d_{x=y}$. For $s\in S,w\in\II_*$ we sometimes set $s\bul w=sw$ if $sw=ws^*$ and $s\bul w=sws^*$ if 
$sw\ne ws^*$; note that $s\bul w\in\II_*$.

For any $s\in S,t\in S,t\ne s$ let $m_{s,t}=m_{t,s}\in[2,\infty]$ be the order of $st$.
For any subset $K$ of $S$ let $W_K$ be the subgroup of $W$ generated by $K$.
If $J\sub K$ are subsets of $S$ we set $W_K^J=\{w\in W_K;l(wy)>l(w)\text{ for any }y\in W_J-\{1\}\}$,
${}^JW_K=\{w\in W_K;l(yw)>l(w)\text{ for any }y\in W_J-\{1\}\}$; note that ${}^JW_K=(W_K^J)\i$.
For any subset $K$ of $S$ such that $W_K$ is finite we denote by $w_K$ the unique element of maximal length of 
$W_K$. 

\head Contents \endhead
1. Involutions and double cosets.

2. Proof of Theorem 0.1.

3. Proof of Theorem 0.2.

4. Proof of Theorem 0.4.

5. The submodule $\uM^K$ of $\uM$.

6. The action of $u\i(T_s+1)$ in the basis $(A_w)$.

7. An inversion formula.

8. A $(-u)$ analogue of weight multiplicities?

9. Reduction modulo $2$.

\head 1. Involutions and double cosets\endhead
\subhead 1.1\endsubhead
Let $K,K'$ be two subsets of $S$ such that $W_K,W_{K'}$ are finite and let $\Om$ be a $(W_K,W_{K'})$-double coset
in $W$. Let $b$ be the unique element of minimal length of $\Om$. Let 
$J=K\cap(bK'b\i)$, $J'=(b\i Kb)\cap K'$ so that $b\i Jb=J'$ hence $b\i W_Jb=W_{J'}$. If $x\in\Om$ then $x=cbd$ 
where $c\in W_K^J$, $d\in W_{K'}$ are uniquely determined; moreover, $l(x)=l(c)+l(b)+l(d)$, see Kilmoyer 
\cite{\KI, Prop. 29}. 
We can write uniquely $d=zc'$ where $z\in W_{J'},c'\in{}^{J'}W_{K'}$; moreover, 
$l(d)=l(z)+l(c')$. Thus we have $x=cbzc'$ where $c\in W_K^J$, $z\in W_{J'},c'\in {}^{J'}W_{K'}$ are uniquely 
determined; moreover, $l(x)=l(c)+l(b)+l(z)+l(c')$. Note that $\tb:=w_Kw_Jbw_{K'}$ is the unique element of 
maximal length of $\Om$; we have $l(\tb)=l(w_K)+l(b)+l(w_{K'})-l(w_J)$.

\subhead 1.2\endsubhead
Now assume in addition that $K'=K^*$ and that $\Om$ is stable under $w\m w^{*-1}$. Then $b^{*-1}\in\Om$, 
$\tb^{*-1}\in\Om$, $l(b^{*-1})=l(b)$, $l(\tb^{*-1})=l(\tb)$, and by uniqueness we have $b^{*-1}=b$, 
$\tb^{*-1}=\tb$, that is, $b\in\II_*$, $\tb\in\II_*$. Also we have $J^*=K^*\cap(b\i Kb)=J'$ hence 
$W_{J'}=(W_J)^*$.
If $x\in\Om\cap\II_*$, then writing $x=cbzc'$ as in 1.1 we 
have $x=x^{*-1}=c'{}^{*-1}b(b\i z^{*-1}b)c^{*-1}$ where 
$c'{}^{*-1}\in({}^JW_{K})\i=W_K^J$, $c^{*-1}\in(W_{K^*}^{J^*})\i={}^{J^*}W_{K^*}$,
$b\i z^{*-1}b\in b\i W_Jb=W_{J^*}$. By the uniqueness of $c,z,c'$, we must have $c'{}^{*-1}=c$, $c^{*-1}=c'$, 
$b\i z^{*-1}b=z$. Conversely, if $c\in W_K^J$, $z\in W_{J^*},c'\in {}^{J^*}W_{K^*}$ are such that $c'{}^{*-1}=c$ 
(hence $c^{*-1}=c'$) and $b\i z^{*-1}b=z$ then clearly $cbzc'\in\Om\cap\II_*$. Note that $y\m b\i y^*b$ is an
automorphism $\t:W_{J^*}@>>>W_{J^*}$ which leaves $J^*$ stable and satisfies $\t^2=1$. Hence 
$\II_\t:=\{y\in W_{J^*};\t(y)\i=y\}$ is well defined. We see that we have a bijection 

(a) $W_K^J\T\II_\t@>>>\Om\cap\II_*$, $(c,z)\m cbzc^{*-1}$.

\subhead 1.3\endsubhead
In the setup of 1.2 we assume that $s\in S$, $K=\{s\}$, so that $K'=\{s^*\}$. In this case we have either

$sb=bs^*$, $J=\{s\}$, $\Om\cap\II_*=\{b,bs^*=\tb\}$, $l(bs^*)=l(b)+1$, or

$sb\ne bs^*$, $J=\em$, $\Om\cap\II_*=\{b,sbs^*=\tb\}$, $l(sbs^*)=l(b)+2$.

\subhead 1.4\endsubhead
In the setup of 1.2 we assume that $s\in S,t\in S,t\ne s$, $m:=m_{s,t}<\infty$,
$K=\{s,t\}$, so that $K^*=\{s^*,t^*\}$. We set $\b=l(b)$. 
For $i\in[1,m]$ we set $\ss_i=sts\do$ ($i$ factors), $\tt_i=tst\do$ ($i$ factors).

We are in one of the following cases (note that we have $sb=bt^*$ if and only if $tb=bs^*$, since $b^{*-1}=b$).

(i) $\{sb,tb\}\cap\{bs^*,bt^*\}=\em$, $J=\em$, 
$\Om\cap\II_*=\{\x_{2i},\x'_{2i} (i\in[0,m]),\x_0=\x'_0=b,\x_{2m}=\x'_{2m}=\ti b\}$
where $\x_{2i}=\ss_i\i b\ss_i^*,\x'_{2i}=\tt_i\i b\tt_i$, $l(\x_{2i})=l(\x'_{2i})=\b+2i$.

(ii) $sb=bs^*$, $tb\ne bt^*$, $J=\{s\}$, $\Om\cap\II_*=\{\x_{2i},\x_{2i+1} (i\in[0,m-1])\}$ where 
$\x_{2i}=\tt_i\i b\tt_i^*$, $l(\x_{2i})=\b+2i$, $\x_{2i+1}=\tt_i\i b\ss_{i+1}^*=\ss_{i+1}\i b\tt_i^*$ , 
$l(\x_{2i+1})=\b+2i+1$, $\x_0=b,\x_{2m-1}=\ti b$.

(iii) $sb\ne bs^*$, $tb=bt^*$, $J=\{t\}$, $\Om\cap\II_*=\{\x_{2i},\x_{2i+1} (i\in[0,m-1])\}$ where 
$\x_{2i}=\ss_i\i b\ss_i^*$, $l(\x_{2i})=\b+2i$, $\x_{2i+1}=\ss_i\i b\tt_{i+1}^*=\tt_{i+1}\i b\ss_i^*$, 
$l(\x_{2i+1})=\b+2i+1$, $\x_0=b,\x_{2m-1}=\ti b$.

(iv) $sb=bs^*$, $tb=bt^*$, $J=K$, $m$ odd, 
$\Om\cap\II_*=\{\x_0=\x'_0=b,\x_{2i+1},\x'_{2i+1} (i\in[0,(m-1)/2]),\x_m=\x'_m=\ti b\}$
where $\x_1=sb$, $\x_3=tstb$, $\x_5=ststsb$, $\do$; $x'_1=tb$, $x'_3=stsb$, $x'_5=tststb$, $\do$; 
$l(\x_{2i+1})=l(\x'_{2i+1})=\b+2i+1$. 

(v) $sb=bs^*$, $tb=bt^*$, $J=K$, $m$ even, 
$\Om\cap\II_*=\{\x_0=\x'_0=b,\x_{2i+1},\x'_{2i+1} (i\in[0,(m-2)/2]),\x_m=\x'_m=\ti b\}$
where $\x_1=sb$, $\x_3=tstb$, $\x_5=ststsb$, $\do$; $\x'_1=tb$, $\x'_3=stsb$, $\x'_5=tststb$, $\do$; 
$l(\x_{2i+1})=l(\x'_{2i+1})=\b+2i+1$, $\x_m=\x'_m=b\ss^*_m=b\tt^*_m=\ss_mb=\tt_mb$, $l(\x_m)=l(\x'_m)=\b+m$.

(vi) $sb=bt^*$, $tb=bs^*$, $J=K$, $m$ odd,
$\Om\cap\II_*=\{\x_0=\x'_0=b,\x_{2i},\x'_{2i} (i\in[0,(m-1)/2]),\x_m=\x'_m=\ti b\}$
where $\x_2=stb$, $\x_4=tstsb$, $\x_6=stststb$, $\do$; $\x'_2=tsb$, $\x'_4=ststb$, $\x'_6=tststsb$, $\do$;
$l(\x_{2i})=l(\x'_{2i})=\b+2i$, $\x_m=\x'_m=b\ss^*_m=b\tt^*_m=\tt_mb=\ss_mb$, $l(\x_m)=l(\x'_m)=\b+m$.

(vii) $sb=bt^*$, $tb=bs^*$, $J=K$, $m$ even,
$\Om\cap\II_*=\{\x_0=\x'_0=b,\x_{2i},\x'_{2i} (i\in[0,m/2]),\x_m=\x'_m=\ti b\}$
where $\x_2=stb$, $\x_4=tstsb$, $\x_6=stststb$, $\do$; $\x'_2=tsb$, $\x'_4=ststb$, $\x'_6=tststsb$, $\do$;
$l(\x_{2i})=l(\x'_{2i})=\b+2i$.

\head 2. Proof of Theorem 0.1\endhead
\subhead 2.1\endsubhead
Let $\dM=\QQ(u)\ot_\ca M$ (a $\QQ(u)$-vector space with basis $\{a_w,w\in\II_*\}$). Let $\dfH=\QQ(u)\ot_\ca\fH$
(a $\QQ(u)$-algebra with basis $\{T_w;w\in W\}$ defined by the relations 0.0(i),(ii)).
The product of a sequence $\x_1,\x_2,\do$ of $k$ elements of $\dfH$ is sometimes denoted by $(\x_1\x_2\do)_k$. 
It is well known that $\dfH$ is the associative $\QQ(u)$-algebra (with $1$) with generators $T_s (s\in S)$ and 
relations 0.0(ii) and

$(T_sT_tT_s\do)_m=(T_tT_sT_t\do)_m$ for any $s\ne t$ in $S$ such that $m:=m_{s,t}<\infty$.
\nl
For $s\in S$ we set $\ocT_s=(u+1)\i(T_s-u)\in\dfH$. Note that $T_s,\ocT_s$ are invertible in $\dfH$: we have
$\ocT_s\i=(u^2-u)\i(T_s+1+u-u^2)$.

\subhead 2.2\endsubhead
For any $s\in S$ we define a $\QQ(u)$-linear map $T_s:\dM@>>>\dM$ by the formulas in 0.1(i)-(iv).
For $s\in S$ we also define a $\QQ(u)$-linear map $\ocT_s:\dM@>>>\dM$ by $\ocT_s=(u+1)\i(T_s-u)$. For $w\in\II_*$
we have: 

(i) $a_{sw}=\ocT_sa_w$ if $sw=ws^*>w$; $a_{sws}=T_sa_w$ if $sw\ne ws^*>w$.

\subhead 2.3\endsubhead
To prove Theorem 0.1 it is enough to show that the formulas 0.1(i)-(iv) define an $\dfH$-module structure on
$\dM$.

Let $s\in S$. To verify that $(T_s+1)(T_s-u^2)=0$ on $\dM$ it is enough to note that the $2\T2$ matrices with 
entries in $\QQ(u)$
$$\left(\sm u&u+1\\u^2-u&u^2-u-1\esm\right)$$   
$$\left(\sm0 &1\\u^2&u^2-1\esm\right)$$   
which represent $T_s$ on the subspace of $\dM$ spanned by $a_w,a_{sw}$ (with $w\in\II,sw=ws^*>w$) or by 
$a_w,a_{sws^*}$ (with $w\in\II,sw\ne ws^*>w$) have eigenvalues $-1,u^2$.

Assume now that $s\ne t$ in $S$ are such that $m:=m_{s,t}<\infty$. It remains to verify the equality
$(T_sT_tT_s\do)_m=(T_tT_sT_t\do)_m:\dM@>>>\dM$. We must show that 
$(T_sT_tT_s\do)_ma_w=(T_tT_sT_t\do)_ma_w$ for any $w\in\II_*$. 
We will do this by reducing the general case to calculations in a dihedral group.

Let $K=\{s,t\}$, so that $K^*=\{s^*,t^*\}$. Let 
$\Om$ be the $(W_K,W_{K^*})$-double coset in $W$ that contains $w$. From the definitions it is clear that the 
subspace $\dM_\Om$ of $\dM$ spanned by $\{a_{w'};w'\in\Om\cap\II_*\}$ is stable under $T_s$ and $T_t$. Hence it 
is enough to show that 

(a) $(T_sT_tT_s\do)_m\mu=(T_tT_sT_t\do)_m\mu$ for any $\mu\in\dM_\Om$. 
\nl 
Since $w^{*-1}=w$ we see that $w'\m w'{}^{*-1}$ maps $\Om$ into itself. Thus $\Om$ is as in 1.2 and we are in 
one of the cases (i)-(vii) in 1.4. The proof of (a) in the various cases is given in 2.4-2.10.
Let $b\in\Om$, $J\sub K$ be as in 1.2. Let $\ss_i,\tt_i$ be as in 1.4.

Let $\dfH_K$ be the subspace of $\dfH$ spanned by $\{T_y;y\in W_K\}$; note that $\dfH_K$ is a $\QQ(u)$-subalgebra
of $\dfH$.

\subhead 2.4\endsubhead
Assume that we are in case 1.4(i). We define an isomorphism of vector spaces $\Ph:\dfH_K@>>>\dM_\Om$ by 
$T_c\m a_{cbc^{*-1}}$ ($c\in W_K$). From definitions we have $T_s\Ph(T_c)=\Ph(T_sT_c)$, $T_t\Ph(T_c)=\Ph(T_tT_c)$
for any $c\in W_K$. It follows that for any $x\in\dfH_K$ we have $T_s\Ph(x)=\Ph(T_sx)$, $T_t\Ph(x)=\Ph(T_tx)$,
hence $(T_sT_tT_s\do)_m\Ph(x)-(T_tT_sT_t\do)_m\Ph(x)=\Ph((T_sT_tT_s\do)_mx-(T_tT_sT_t\do)_mx)=0$.
(We use that $(T_sT_tT_s\do)_m=(T_tT_sT_t\do)_m$ in $\dfH_K$.) Since $\Ph$ is an isomorphism we deduce that 
2.3(a) holds in our case.

Assume that we are in case 1.4(ii). We define $r,r'$ by $r=s$, $r'=t$ if $m$ is odd, $r=t$, $r'=s$ if $m$ is even.
We have 

$a_{\x_0}@>T_t>>a_{\x_2}@>T_s>>a_{\x_4}@>T_t>>\do@>T_r>>a_{\x_{2m-2}}$,

$a_{\x_1}@>T_t>>a_{\x_3}@>T_s>>a_{\x_5}@>T_s>>\do@>T_r>>a_{\x_{2m-1}}$.
\nl
We have $s\x_0=\x_0s^*=\x_1$ hence $a_{\x_0}@>T_s>>ua_{\x_0}+(u+1)a_{\x_1}$. We show that

$r'\x_{2m-2}=\x_{2m-2}r'{}^*=\x_{2m-1}$ 
\nl
We have $r'\x_{2m-2}=\do tstbt^*s^*t^*\do$ where the product to the left (resp. right) of $b$ has $m$ (resp. 
$m-1$) factors). Using the definition of $m$ and the identity $sb=bs^*$ we deduce 
$r'\x_{2m-2}=\do stsbt^*s^*t^*\do=\do stbs^*t^*s^*\do$ (in the last expression the product to the left (resp. 
right) of $b$ has $m-1$ (resp. $m$) factors). Thus $r'\x_{2m-2}=\x_{2m-1}$. Using again the definition of $m$ we 
have $\x_{2m-1}=\do stbt^*s^*t^*\do$ where the product to the left (resp. right) of $b$ has $m-1$ (resp. $m$) 
factors. Thus  $\x_{2m-1}=\x_{2m-2}r'{}^*$ as required.

We deduce that

$a_{\x_{2m-2}}@>T_{r'}>>ua_{\x_{2m-2}}+(u+1)a_{\x_{2m-1}}$.
\nl
We set $a'_{\x_1}=ua_{\x_0}+(u+1)a_{\x_1}$, $a'_{\x_3}=ua_{\x_2}+(u+1)a_{\x_3}$, $\do$,
$a'_{\x_{2m-1}}=ua_{\x_{2m-2}}+(u+1)a_{\x_{2m-1}}$. Note that $a_{\x_0},a_{\x_2},a_{\x_4},\do,a_{\x_{2m-2}}$ 
together with $a'_{\x_1},a'_{\x_3},\do,a'_{\x_{2m-1}}$ form a basis of $\dM_\Om$ and we have 

$a_{\x_0}@>T_t>>a_{\x_2}@>T_s>>a_{\x_4}@>T_t>>\do@>T_r>>a_{\x_{2m-2}}@>T_{r'}>>a'_{\x_{2m-1}}$

$a_{\x_0}@>T_s>>a'_{\x_1}@>T_t>>a'_{\x_3}@>T_s>>a'_{\x_5}@>T_s>>\do @>T_r>>a'_{\x_{2m-1}}$.
\nl
We define an isomorphism of vector spaces $\Ph:\dfH_K@>>>\dM_\Om$ by 
$1\m a_{\x_0}$, $T_t\m a_{\x_2}$, $T_sT_t\m a_{\x_4}$, $\do$, $T_r\do T_sT_t\m a_{\x_{2m-2}}$ (the product has
$m-1$ factors), $T_s\m\a'_{\x_1}$, $T_tT_s\m a'_{\x_3}$, $\do$, $T_r\do T_tT_s\m a'_{\x_{2m-1}}$ (the product has
$m$ factors). From definitions for any $c\in W_K$ we have 

(a) $T_s\Ph(T_c)=\Ph(T_sT_c)$ if $sc>c$, $T_t\Ph(T_c)=\Ph(T_tT_c)$ if $tc>c$,

(b) $T_s\i\Ph(T_c)=\Ph(T_s\i T_c)$ if $sc<c$, $T_t\i\Ph(T_c)=\Ph(T_t\i T_c)$ if $tc<c$.
\nl
Since $T_s=u^2T_s\i+(u^2-1)$ both as endomorphisms of $\dM$ and as elements of $\dfH$ we see that (b) implies 
that $T_s\Ph(T_c)=\Ph(T_sT_c)$ if $sc<c$. Thus $T_s\Ph(T_c)=\Ph(T_sT_c)$ for any $c\in W_K$. Similarly,
$T_t\Ph(T_c)=\Ph(T_tT_c)$ for any $c\in W_K$. It follows that for any $x\in\dfH_K$ we have $T_s\Ph(x)=\Ph(T_sx)$, 
$T_t\Ph(x)=\Ph(T_tx)$, hence $(T_sT_tT_s\do)\Ph(x)-(T_tT_sT_t\do)\Ph(x)=\Ph((T_sT_tT_s\do)x-(T_tT_sT_t\do)x)=0$ 
where the products $T_sT_tT_s\do,T_tT_sT_t\do$ have $m$ factors. (We use that $T_sT_tT_s\do=T_tT_sT_t\do$ in 
$\dfH_K$.) Since $\Ph$ is an isomorphism we deduce that $(T_sT_tT_s\do)\mu-(T_tT_sT_t\do)\mu=0$ for any 
$\mu\in\dM_\Om$. Hence 2.3(a) holds in our case.

\subhead 2.5\endsubhead
Assume that we are in case 1.4(iii). By the argument in case 1.4(ii) with $s,t$ interchanged we see that (a) 
holds in our case.

\subhead 2.6\endsubhead
Assume that we are in one of the cases 1.4(iv)-(vii). We have $J=K$ that is, $K=bK^*b\i$. We have
$\Om=W_Kb=bW_{K^*}$. Define $m'\ge1$ by $m=2m'+1$ if $m$ is odd, $m=2m'$ if $m$ is even.
Define $s',t'$ by $s'=s,t'=t$ if $m'$ is even, $s'=t,t'=s$ if $m'$ is odd. 

\subhead 2.7\endsubhead
Assume that we are in case 1.4(iv). We define some elements of $\dfH_K$ as follows:
$$\align&\et_0=T_{\ss_{m'}}+T_{\tt_{m'}}+(1+u-u^2)(T_{\ss_{m'-1}}+T_{\tt_{m'-1}})\\&
+(1+u-u^2-u^3+u^4)(T_{\ss_{m'-2}}+T_{\tt_{m'-2}})+\do\\&
+(1+u-u^2-u^3+u^4+u^5-\do+(-1)^{m'-2}u^{2m'-4}+(-1)^{m'-2}u^{2m'-3}\\&+(-1)^{m'-1}u^{2m'-2})(T_{\ss_1}+T_{\tt_1})
\\&+(1+u-u^2-u^3+u^4+u^5-\do+(-1)^{m'-1}u^{2m'-2}\\&+(-1)^{m'-1}u^{2m'-1}+(-1)^{m'}u^{2m'}),\endalign$$
$$\et_1=\ocT_s\et_0, \et_3=T_t\et_1,\do, \et_{2m'-1}=T_{t'}\et_{2m'-3}, \et_{2m'+1}=T_{s'}\et_{2m'-1},$$
$$\et'_1=\ocT_t\et_0, \et'_3=T_s\et'_1,\do, \et'_{2m'-1}=T_{s'}\et'_{2m'-3},\et'_{2m'+1}=T_{t'}\et'_{2m'-1}.$$
For example if $m=7$ we have 
$$\align&\et_0=T_{sts}+T_{tst}+(1+u-u^2)T_{ts}+(1+u-u^2)T_{st}+(1+u-u^2-u^3+u^4)T_s\\&+(1+u-u^2-u^3+u^4)T_t+
(1+u-u^2-u^3+u^4+u^5-u^6),\endalign$$
$$\align&\et_1=(u+1)\i(T_{stst}-uT_{tst}+(-u+u^3)T_{ts}+(-u+u^3)T_{st}
+(-u+2u^3-u^5)T_s\\&+(-u+2u^3-u^5)T_t+(-u+2u^3-2u^5+u^7)),\endalign$$
$$\align&\et'_1=(u+1)\i(T_{tsts}-uT_{sts}+(-u+u^3)T_{ts}+(-u+u^3)T_{st}
+(-u+2u^3-u^5)T_s\\&+(-u+2u^3-u^5)T_t+(-u+2u^3-2u^5+u^7)),\endalign$$
$$\et_3=(u+1)\i(T_{tstst}-u^3T_{st}+(-u^3+u^5)T_s+(-u^3+u^5)T_t+(-u^3+2u^5-u^7)),$$
$$\et'_3=(u+1)\i(T_{ststs}-u^3T_{ts}+(-u^3+u^5)T_s+(-u^3+u^5)T_t+(-u^3+2u^5-u^7)),$$
$$\et_5=(u+1)\i(T_{ststst}-u^5T_t+(-u^5+u^7)),$$
$$\et'_5=(u+1)\i(T_{tststs}-u^5T_s+(-u^5+u^7)),$$
$$\et_7=\et'_7=(u+1)\i(T_{stststs}-u^7).$$     
One checks by direct computation in $\dfH_K$ that
$$\et_m=\et'_m=(u+1)\i(T_{\ss_m}-u^m)\tag a$$
and that the elements $\et_0,\et_1,\et'_1,\et_3,\et'_3,\do\et_{2m'-1},\et'_{2m'-1},\et_m$
are linearly independent in $\dfH_K$; they span a subspace of $\dfH_K$ denoted by $\dfH_K^+$. From (a) we deduce:
$$(T_{s'}T_{t'}T_{s'}\do T_tT_sT_t\ocT_s)_{m'+1}\et_0=(T_{t'}T_{s'}T_{t'}\do T_sT_tT_s\ocT_t)_{m'+1}\et_0.\tag b$$
We have
$$\ocT_s\i\et_1=\et_0,T_t\i\et_3=\et_1,\do,T_{t'}\i\et_{2m'-1}=\et_{2m'-3},T_{s'}\i\et_{2m'+1}=\et_{2m'-1},$$
$$\ocT_t\i\et'_1=\et_0,T_s\i\et'_3=\et'_1,\do,T_{s'}\i\et'_{2m'-1}=\et'_{2m'-3},T_{t'}\i\et'_{2m'+1}
=\et'_{2m'-1}.$$
It follows that $\dfH_K^+$ is stable under left multiplication by $T_s$ and $T_t$ hence it is a left ideal of
$\dfH_K$.
From the definitions we have
$$a_{\x_1}=\ocT_sa_{\x_0}, a_{\x_3}=T_ta_{\x_1},\do, a_{\x_{2m'-1}}=T_{t'}a_{\x_{2m'-3}}, 
a_{\x_{2m'+1}}=T_{s'}a_{\x_{2m'-1}},$$
$$a_{\x'_1}=\ocT_ta_{\x_0}, a_{\x'_3}=T_sa_{\x'_1},\do, a_{\x'_{2m'-1}}=T_{s'}a_{\x'_{2m'-3}},
a_{\x'_{2m'+1}}=T_{t'}a_{\x'_{2m'-1}},$$
$$\ocT_s\i a_{\x_1}=a_{\x_0},T_t\i a_{\x_3}=a_{\x_1},\do,T_{t'}\i a_{\x_{2m'-1}}=a_{\x_{2m'-3}},
T_{s'}\i a_{\x_{2m'+1}}=a_{\x_{2m'-1}},$$
$$\ocT_t\i a_{\x'_1}=a_{\x_0},T_s\i a_{\x'_3}=a_{\x'_1},\do,T_{s'}\i a_{\x'_{2m'-1}}=a_{\x'_{2m'-3}},
T_{t'}\i a_{\x'_{2m'+1}}=a_{\x'_{2m'-1}}.$$
Hence the vector space isomorphism $\Ph:\dfH_K^+@>\si>>\dM_\Om$ given by $\et_{2i+1}\m a_{\x_{2i+1}}$,
$\et'_{2i+1}\m a_{\x'_{2i+1}}$ $(i\in[0,(m-1)/2])$, $\et_0\m a_{\x_0}$ satisfies $\Ph(T_sh)=T_s\Ph(h)$,
$\Ph(T_th)=T_t\Ph(h)$ for any $h\in\dfH_K^+$. 
Since $(T_sT_tT_s\do)_mh=(T_tT_sT_t\do)_mh$ for $h\in\dfH_K^+$, we deduce that  2.3(a) holds in our case.

\subhead 2.8\endsubhead
Assume that we are in case 1.4(v). We define some elements of $\dfH_K$ as follows:
$$\align&\et_0=T_{\ss_{m'-1}}+T_{\tt_{m'-1}}+(1-u^2)(T_{\ss_{m'-2}}+T_{\tt_{m'-2}})\\&
+(1-u^2+u^4)(T_{\ss_{m'-3}}+T_{\tt_{m'-3}})+\do\\&
+(1-u^2+u^4-\do+(-1)^{m'-2}u^{2(m'-2)})(T_{\ss_1}+T_{\tt_1})\\&
+(1-u^2+u^4-\do+(-1)^{m'-1}u^{2(m'-1)}),\endalign$$
(if $m\ge4$), $\et_0=1$ (if $m=2$),
$$\et_1=\ocT_s\et_0, \et_3=T_t\et_1,\do, \et_{2m'-1}=T_{t'}\et_{2m'-3}, \et_{2m'}=\ocT_{s'}\et_{2m'-1},$$
$$\et'_1=\ocT_t\et_0, \et'_3=T_s\et'_1,\do, \et'_{2m'-1}=T_{s'}\et'_{2m'-3}, \et'_{2m'}=\ocT_{t'}\et'_{2m'-1}.$$
For example if $m=4$ we have
$$\et_0=T_s+T_t+(1-u^2),$$
$$\et_1=(u+1)\i(T_{st}-uT_s-uT_t+(-u+u^2+u^3)),$$  
$$\et'_1=(u+1)\i(T_{ts}-uT_s-uT_t+(-u+u^2+u^3)),$$
$$\et_3=(u+1)\i(T_{tst}-uT_{ts}+u^2T_s-u^3),$$
$$\et'_3=(u+1)\i(T_{sts}-uT_{st}+u^2T_s-u^3),$$
$$\et_4=\et'_4=(u+1)^{-2}(T_{stst}-uT_{sts}-uT_{tst}+u^2T_{st}+u^2T_{ts}-u^3T_s-u^3T_t+u^4).$$
If $m=6$ we have
$$\et_0=T_{st}+T_{ts}+(1-u^2)T_s+(1-u^2)T_t+(1-u^2+u^4),$$
$$\align&\et_1=(u+1)\i(T_{sts}-uT_{st}-uT_{ts}+(-u+u^2+u^3)T_s\\&+(-u+u^2+u^3)T_t+(-u+u^2+u^3-u^4-u^5)),\endalign
$$
$$\align&\et'_1=(u+1)\i(T_{tst}-uT_{st}-uT_{ts}+(-u+u^2+u^3)T_s\\&+(-u+u^2+u^3)T_t+(-u+u^2+u^3-u^4-u^5)),\endalign
$$
$$\et_3=(u+1)\i(T_{tsts}-uT_{tst}-u^2T_{ts}-u^3T_s-u^3T_t+(-u^3+u^4+u^5)),$$
$$\et'_3=(u+1)\i(T_{stst}-uT_{sts}-u^2T_{st}-u^3T_s-u^3T_t+(-u^3+u^4+u^5)),$$
$$\et_5=(u+1)\i(T_{ststs}-uT_{stst}-u^2T_{sts}-u^3T_{st}+u^4T_s-u^5),$$
$$\et'_5=(u+1)\i(T_{tstst}-uT_{tsts}-u^2T_{tst}-u^3T_{ts}+u^4T_t-u^5),$$
$$\align&\et_6=\et'_6=(u+1)^{-2}(T_{ststst}-uT_{ststs}-uT_{tstst}+u^2T_{stst}+u^2T_{tsts}-u^3T_{sts}\\&
-u^3T_{tst}+u^4T_{st}+u^4T_{ts}-u^5T_s-u^5T_t+u^6).\endalign$$
If $m=8$ we have
$$\align&\et_0=T_{sts}+T_{tst}+(1-u^2)T_{st}+(1-u^2)T_{ts}+(1-u^2+u^4)T_s+(1-u^2+u^4)T_t\\&
+(1-u^2+u^4-u^6),\endalign$$
$$\align&\et_1=(u+1)\i(T_{stst}-uT_{sts}-uT_{tst}+(-u+u^2+u^3)T_{st}+(-u+u^2+u^3)T_{ts}\\&+(-u+u^2+u^3-u^4-u^5)T_s
+(-u+u^2+u^3-u^4-u^5)T_t\\&+(-u+u^2+u^3-u^4-u^5+u^6+u^7)),\endalign$$
$$\align&\et'_1=(u+1)\i(T_{tsts}-uT_{sts}-uT_{tst}+(-u+u^2+u^3)T_{st}+(-u+u^2+u^3)T_{ts}\\&
+(-u+u^2+u^3-u^4-u^5)T_s+(-u+u^2+u^3-u^4-u^5)T_t+\\&(-u+u^2+u^3-u^4-u^5+u^6+u^7)),\endalign$$
$$\align&\et_3=(u+1)\i(T_{tstst}-uT_{tsts}+u^2T_{tst}-u^3T_{st}-u^3T_{ts}+(-u^3+u^4+u^5)T_s\\&+(-u^3+u^4+u^5)T_t+
(-u^3+u^4+u^5-u^6-u^7)),\endalign$$
$$\align&\et'_3=(u+1)\i(T_{ststs}-uT_{stst}+u^2T_{sts}-u^3T_{st}-u^3T_{ts}+(-u^3+u^4+u^5)T_s\\&+(-u^3+u^4+u^5)T_t+
(-u^3+u^4+u^5-u^6-u^7)),\endalign$$
$$\align&\et_5=(u+1)\i(T_{ststst}-uT_{ststs}+u^2T_{stst}-u^3T_{sts}+u^4T_{st}-u^5T_s-u^5T_t\\&
+(-u^5+u^6+u^7)),\endalign$$
$$\align&\et'_5=(u+1)\i(T_{tststs}-uT_{tstst}+u^2T_{tsts}-u^3T_{tst}+u^4T_{ts}-u^5T_s-u^5T_t\\&
+(-u^5+u^6+u^7)),\endalign$$
$$\align&\et_7=(u+1)\i(T_{tststst}-uT_{tststs}+u^2T_{tstst}-u^3T_{tsts}+u^4T_{tst}-u^5T_{ts}\\&+u^6T_t-u^7),
\endalign$$
$$\align&\et'_7=(u+1)\i(T_{stststs}-uT_{ststst}+u^2T_{ststs}-u^3T_{stst}+u^4T_{sts}-u^5T_{st}\\&+u^6T_s-u^7),
\endalign$$
$$\align&\et_8=\et'_8=(u+1)^{-2}(T_{stststst}-uT_{stststs}-uT_{tststst}+u^2T_{tststs}\\&+u^2T_{ststst}
-u^3T_{ststs}-u^3T_{tstst}+u^4T_{stst}+u^4T_{tsts}-u^5T_{sts}-u^5T_{tst}+u^6T_{st}\\&
+u^6T_{ts}-u^7T_s-u^tT_t+u^8).\endalign$$
One checks by direct computation in $\dfH_K$ that
$$\et_m=\et'_m=(u+1)^{-2}\sum_{y\in W_K}(-u)^{m-l(y)}T_y\tag a$$
and that the elements $\et_0,\et_1,\et'_1,\et_3,\et'_3,\do\et_{2m'-1},\et'_{2m'-1},\et_m$
are linearly independent in $\dfH_K$; they span a subspace of $\dfH_K$ denoted by $\dfH_K^+$. From (a) we deduce:
$$(\ocT_{s'}T_{t'}T_{s'}\do T_tT_sT_t\ocT_s)_{m'+1}\et_0=(\ocT_{t'}T_{s'}T_{t'}\do T_sT_tT_s\ocT_t)_{m'+1}\et_0.
\tag b$$
We have
$$\ocT_s\i\et_1=\et_0,T_t\i\et_3=\et_1,\do,T_{t'}\i\et_{2m'-1}=\et_{2m'-3},\ocT_{s'}\i\et_{2m'}=\et_{2m'-1},$$
$$\ocT_t\i\et'_1=\et_0,T_s\i\et'_3=\et'_1,\do,T_{s'}\i\et'_{2m'-1}=\et'_{2m'-3},\ocT_{t'}\i\et'_{2m'}
=\et'_{2m'-1}.$$
It follows that $\dfH_K^+$ is stable under left multiplication by $T_s$ and $T_t$ hence it is a left ideal of
$\dfH_K$.
From the definitions we have
$$a_{\x_1}=\ocT_sa_{\x_0}, a_{\x_3}=T_ta_{\x_1},\do, a_{\x_{2m'-1}}=T_{t'}a_{\x_{2m'-3}}, 
a_{\x_{2m'}}=\ocT_{s'}a_{\x_{2m'-1}},$$
$$a_{\x'_1}=\ocT_ta_{\x_0}, a_{\x'_3}=T_sa_{\x'_1},\do, a_{\x'_{2m'-1}}=T_{s'}a_{\x'_{2m'-3}},
a_{\x'_{2m'}}=\ocT_{t'}a_{\x'_{2m'-1}},$$
$$\ocT_s\i a_{\x_1}=a_{\x_0},T_t\i a_{\x_3}=a_{\x_1},\do,T_{t'}\i a_{\x_{2m'-1}}=a_{\x_{2m'-3}},
\ocT_{s'}\i a_{\x_{2m'}}=a_{\x_{2m'-1}},$$
$$\ocT_t\i a_{\x'_1}=a_{\x_0},T_s\i a_{\x'_3}=a_{\x'_1},\do,T_{s'}\i a_{\x'_{2m'-1}}=a_{\x'_{2m'-3}},
\ocT_{t'}\i a_{\x'_{2m'}}=a_{\x'_{2m'-1}}.$$
Hence the vector space isomorphism $\Ph:\dfH_K^+@>\si>>\dM_\Om$ given by $\et_{2i+1}\m a_{\x_{2i+1}}$,
$\et'_{2i+1}\m a_{\x'_{2i+1}}$ $(i\in[0,(m-2)/2])$, $\et_0\m a_{\x_0}$, $\et_m\m a_{\x_m}$ satisfies 
$\Ph(T_sh)=T_s\Ph(h)$, $\Ph(T_th)=T_t\Ph(h)$ for any $h\in\dfH_K^+$. 
Since $(T_sT_tT_s\do)_mh=(T_tT_sT_t\do)_mh$ for $h\in\dfH_K^+$, we deduce that 2.3(a) holds in our case.

\subhead 2.9\endsubhead
Assume that we are in case 1.4(vi). We define some elements of $\dfH_K$ as follows:
$$\align&\et_0=T_{\ss_{m'}}+T_{\tt_{m'}}+(1-u-u^2)(T_{\ss_{m'-1}}+T_{\tt_{m'-1}})\\&
+(1-u-u^2+u^3+u^4)(T_{\ss_{m'-2}}+T_{\tt_{m'-2}})+\do\\&
+(1-u-u^2+u^3+u^4-u^5-\do+(-1)^{m'-2}u^{2m'-4}+(-1)^{m'-1}u^{2m'-3}\\&
+(-1)^{m'-1}u^{2m'-2})(T_{\ss_1}+T_{\tt_1})\\&
+(1+u-u^2-u^3+u^4+u^5-\do+(-1)^{m'-1}u^{2m'-2}\\&+(-1)^{m'}u^{2m'-1}+(-1)^{m'}u^{2m'}),\endalign$$
$$\et_2=T_s\et_0, \et_4=T_t\et_2,\do, \et_{2m'}=T_{s'}\et_{2m'-2}, \et_{2m'+1}=\ocT_{t'}\et_{2m'},$$
$$\et'_2=T_t\et_0, \et'_4=T_s\et'_2,\do, \et'_{2m'}=T_{t'}\et'_{2m'-2}, \et'_{2m'+1}=\ocT_{s'}\et'_{2m'}.$$
For example if $m=7$ we have
$$\align&\et_0=T_{sts}+T_{tst}+(1-u-u^2)T_{ts}+(1-u-u^2)T_{st}+(1-u-u^2+u^3+u^4)T_s\\&+(1-u-u^2+u^3+u^4)T_t+
(1-u-u^2+u^3+u^4-u^5-u^6),\endalign$$
$$\et_2=T_{stst}-uT_{sts}+u^2T_{st}+(u^2-u^3-u^4)T_s+(u^2-u^3-u^4)T_t+(u^2-u^3-u^4+u^5+u^6),$$
$$\et'_2=T_{tsts}-uT_{tst}+u^2T_{st}+(u^2-u^3-u^4)T_s+(u^2-u^3-u^4)T_t+(u^2-u^3-u^4+u^5+u^6),$$
$$\et_4=T_{ststs}-uT_{stst}+u^2T_{sts}-u^3T_{st}+u^4T_s+u^4T_t+(u^4-u^5-u^6),$$
$$\et'_4=T_{tstst}-uT_{tsts}+u^2T_{tst}-u^3T_{ts}+u^4T_s+u^4T_t+(u^4-u^5-u^6),$$
$$\et_6=T_{ststst}-uT_{ststs}+u^2T_{stst}-u^3T_{sts}+u^4T_{st}-u^5T_s+u^6,$$
$$\et'_6=T_{tststs}-uT_{tstst}+u^2T_{tsts}-u^3T_{tst}+u^4T_{ts}-u^5T_t+u^6,$$
$$\align&
\et_7=\et'_7=(u+1)\i(T_{stststs}-uT_{ststst}-uT_{tststs}+u^2T_{ststs}+u^2T_{tstst}-u^3T_{stst}\\&-u^3T_{tsts}
+u^4T_{sts}+u^4T_{tst}-u^5T_{st}-u^5T_{ts}+u^6T_s+u^6T_t-u^7).\endalign$$ 
One checks by direct computation in $\dfH_K$ that
$$\et_m=\et'_m=(u+1)\i\sum_{y\in W_K}(-u)^{m-l(y)}T_y\tag a$$
and that the elements $\et_0,\et_2,\et'_2,\et_4,\et'_4,\do\et_{2m'},\et'_{2m'},\et_m$
are linearly independent in $\dfH_K$; they span a subspace of $\dfH_K$ denoted by $\dfH_K^+$. From (a) we deduce:
$$(\ocT_{s'}T_{t'}T_{s'}\do T_tT_s)_{m'+1}\et_0=(\ocT_{t'}T_{s'}T_{t'}\do T_sT_t)_{m'+1}\et_0.\tag b$$
We have
$$\et_0=T_s\i\et_2, \et_2=T_t\i\et_4,\do, \et_{2m'-2}=T_{s'}\i\et_{2m'}, \et_{2m'}=\ocT_{t'}\i\et_{2m'+1},$$
$$\et_0=T_t\i\et'_2, \et'_2=T_s\i\et'_4,\do,\et'_{2m'-2}=T_{t'}\i\et'_{2m'}, \et'_{2m'}=\ocT_{s'}\i\et_{2m'+1}.$$
It follows that $\dfH_K^+$ is stable under left multiplication by $T_s$ and $T_t$ hence it is a left ideal of
$\dfH_K$.
From the definitions we have
$$a_{\x_2}=T_sa_{\x_0}, a_{\x_4}=T_ta_{\x_2},\do, a_{\x_{2m'}}=T_{s'}a_{\x_{2m'-2}},
a_{\x_{2m'+1}}=\ocT_{t'}a_{\x_{2m'}},$$
$$a_{\x'_2}=T_ta_{\x_0}, a_{\x'_4}=T_sa_{\x'_2},\do, a_{\x'_{2m'}}=T_{t'}a_{\x'_{2m'-2}}, 
a_{\x'_{2m'+1}}=\ocT_{s'}a_{\x_{2m'}},$$
$$a_{\x_0}=T_s\i a_{\x_2}, a_{\x_2}=T_t\i a_{\x_4},\do, a_{\x_{2m'-2}}=T_{s'}\i a_{\x_{2m'}},
a_{\x_{2m'}}=\ocT_{t'}\i a_{\x_{2m'+1}},$$
$$a_{\x_0}=T_t\i a_{\x'_2}, a_{\x'_2}=T_s\i a_{\x'_4},\do, a_{\x'_{2m'-2}}=T_{t'}\i a_{\x'_{2m'}}, 
a_{\x'_{2m'}}=\ocT_{s'}\i a_{\x_{2m'+1}}.$$
Hence the vector space isomorphism $\Ph:\dfH_K^+@>\si>>\dM_\Om$ given by $\et_{2i}\m a_{\x_{2i}}$,
$\et'_{2i}\m a_{\x'_{2i}}$ $(i\in[0,(m-1)/2])$, $\et_m\m a_{\x_m}$ satisfies 
$\Ph(T_sh)=T_s\Ph(h)$, $\Ph(T_th)=T_t\Ph(h)$ for any $h\in\dfH_K^+$. 
Since $(T_sT_tT_s\do)_mh=(T_tT_sT_t\do)_mh$ for $h\in\dfH_K^+$, we deduce that  2.3(a) holds in our case.

\subhead 2.10\endsubhead
Assume that we are in case 1.4(vii). We define some elements of $\dfH_K$ as follows:
$$\align&\et_0=T_{\ss_{m'}}+T_{\tt_{m'}}+(1-u^2)(T_{\ss_{m'-1}}+T_{\tt_{m'-1}})\\&
+(1-2u^2+u^4)(T_{\ss_{m'-3}}+T_{\tt_{m'-3}})+\do\\&
+(1-2u^2+2u^4-\do+(-1)^{m'-2}2u^{2(m'-2)}+(-1)^{m'-1}u^{2(m'-1)})(T_{\ss_1}+T_{\tt_1})\\&+
(1-2u^2+2u^4-\do+(-1)^{m'-1}2u^{2(m'-1)}+(-1)^{m'}u^{2m'}),\endalign$$
$$\et_2=T_s\et_0, \et_4=T_t\et_2,\do, \et_{2m'}=T_{s'}\et_{2m'-2},$$
$$\et'_2=T_t\et_0, \et'_4=T_s\et'_2,\do, \et'_{2m'}=T_{t'}\et'_{2m'-2}.$$
For example if $m=8$ we have
$$\align&\et_0=T_{stst}+T_{tsts}+(1-u^2)T_{sts}+(1-u^2)T_{tst}+(1-2u^2+u^4)T_{st}\\&+(1-2u^2+u^4)T_{ts}
+(1-2u^2+2u^4-u^6)T_s+(1-2u^2+2u^4-u^6)T_t+\\&(1-2u^2+2u^4-2u^6+u^8),\endalign$$
$$\align&\et_2=T_{ststs}+u^2T_{tst}+(u^2-u^4)T_{st}+(u^2-u^4)T_{ts}+(u^2-2u^4+u^6)T_s\\&+(u^2-2u^4+u^6)T_t+
(u^2-2u^4+2u^6-u^8),\endalign$$
$$\align&\et'_2=T_{tstst}+u^2T_{sts}+(u^2-u^4)T_{st}+(u^2-u^4)T_{ts}+(u^2-2u^4+u^6)T_s\\&+(u^2-2u^4+u^6)T_t+
(u^2-2u^4+2u^6-u^8),\endalign$$
$$\et_4=T_{tststs}+u^4T_{st}+(u^4-u^6)T_s+(u^4-u^6)T_t+(u^4-2u^6+u^8),$$
$$\et'_4=T_{ststst}+u^4T_{ts}+(u^4-u^6)T_s+(u^4-u^6)T_t+(u^4-2u^6+u^8),$$
$$\et_6=T_{stststs}+u^6T_t+(u^6-u^8),$$
$$\et'_6=T_{tststst}+u^6T_s+(u^6-u^8),$$
$$\et_8=\et'_8=T_{stststst}+u^8.$$
One checks by direct computation in $\dfH_K$ that
$$\et_m=\et'_m=T_{\ss_m}+u^m\tag a$$
and that the elements $\et_0,\et_2,\et'_2,\et_4,\et'_4,\do\et_{2m'},\et'_{2m'},\et_m$
are linearly independent in $\dfH_K$; they span a subspace of $\dfH_K$ denoted by $\dfH_K^+$. From (a) we deduce:
$$(T_{t'}T_{s'}\do T_tT_s)_{m'}\et_0=(T_{s'}T_{t'}\do T_sT_t)_{m'}\et_0.\tag c$$
We have
$$\et_0=T_s\i\et_2, \et_2=T_t\i\et_4,\do, \et_{2m'-2}=T_{s'}\i\et_{2m'},$$
$$\et_0=T_t\i\et'_2, \et'_2=T_s\i\et'_4,\do, \et'_{2m'-2}=T_{t'}\i\et'_{2m'}.$$
It follows that $\dfH_K^+$ is stable under left multiplication by $T_s$ and $T_t$ hence it is a left ideal of
$\dfH_K$.
From the definitions we have
$$a_{\x_2}=T_sa_{\x_0}, a_{\x_4}=T_ta_{\x_2},\do, a_{\x_{2m'}}=T_{s'}a_{\x_{2m'-2}},$$
$$a_{\x'_2}=T_ta_{\x_0}, a_{\x'_4}=T_sa_{\x'_2},\do, a_{\x'_{2m'}}=T_{t'}a_{\x'_{2m'-2}},$$
$$a_{\x_0}=T_s\i a_{\x_2}, a_{\x_2}=T_t\i a_{\x_4},\do, a_{\x_{2m'-2}}=T_{s'}\i a_{\x_{2m'}},$$
$$a_{\x_0}=T_t\i a_{\x'_2}, a_{\x'_2}=T_s\i a_{\x'_4},\do, a_{\x'_{2m'-2}}=T_{t'}\i a_{\x'_{2m'}}.$$
Hence the vector space isomorphism $\Ph:\dfH_K^+@>\si>>\dM_\Om$ given by $\et_{2i}\m a_{\x_{2i}}$,
$\et'_{2i}\m a_{\x'_{2i}}$ $(i\in[0,m/2])$ satisfies 
$\Ph(T_sh)=T_s\Ph(h)$, $\Ph(T_th)=T_t\Ph(h)$ for any $h\in\dfH_K^+$. Since $(T_sT_tT_s\do)_mh=(T_tT_sT_t\do)_mh$ 
for $h\in\dfH_K^+$, we deduce that 2.3(a) holds in our case.
This completes the proof of Theorem 0.1.

\subhead 2.11\endsubhead
We show that the $\dfH$-module 
$\dM$ is generated by $a_1$. Indeed, from 2.2(i) we see by induction on $l(w)$ that 
for any $w\in\II_*$, $a_w$ belongs to the $\dfH$-submodule of $\dM$ generated by $a_1$.

\head 3. Proof of Theorem 0.2\endhead
\subhead 3.1\endsubhead
We define a $\ZZ$-linear map $B:M@>>>M$ by $B(u^na_w)=\e_wu^{-n}T_{w^*}\i a_{w^*}$ for any $w\in\II_*,n\in\ZZ$. 
Note that $B(a_1)=a_1$. 

For any $w\in\II_*,s\in S$ we show:

(a) $B(T_sa_w)=T_s\i B(a_w)$.
\nl
Assume first that $sw=ws^*>w$. We must show that $B(ua_w+(u+1)a_{sw})=T_s\i B(a_w)$ or that
$$u\i\e_wT_{w^*}\i a_{w^*}-(u\i+1)\e_wT_{s^*w^*}\i a_{s^*w^*}=T_s\i\e_wT_{w^*}\i a_{w^*}$$
or that
$$T_{w^*}\i a_{w^*}-(u+1)T_{w^*}\i T_{s^*}\i a_{s^*w^*}=u  T_{w^*}\i T_{s^*}\i a_{w^*}$$
or that
$$T_{s^*}a_{w^*}-(u+1) a_{s^*w^*}=ua_{w^*}.$$
This follows from 0.1(i) with $s,w$ replaced by $s^*,w^*$.

Assume next $sw=ws^*<w$. We set $y=sw\in\II_*$ so that $sy>y$. We must show that 
$B((u^2-u-1)a_{sy}+(u^2-u)a_y)=T_s\i B(a_{sy})$ or that
$$-(u^{-2}-u\i-1)\e_yT_{s^*y^*}\i a_{s^*y^*}+(u^{-2}-u\i)\e_yT_{y^*}\i a_{y^*}=-T_s\i\e_yT_{s^*y^*}\i a_{s^*y^*}$$
or that
$$-(u^{-2}-u\i-1)T_{y^*}\i T_{s^*}\i a_{s^*y^*}+(u^{-2}-u\i)T_{y^*}\i a_{y^*}=-T_{y^*}\i T_{s^*}^{-2}a_{s^*y^*}$$ 
or that 
$$-(u^{-2}-u\i-1)T_{s^*}\i a_{s^*y^*}+(u^{-2}-u\i)a_{y^*}=-T_{s^*}^{-2}a_{s^*y^*}$$
or that 
$$-(1-u-u^2)a_{s^*y^*}+(1-u)T_{s^*}a_{y^*}=-(T_{s^*}+1-u^2) a_{s^*y^*}.$$
Using 0.1(i),(ii) with $w,s$ replaced by $y^*,s^*$ we see that it is enough to show that
$$\align&-(1-u-u^2)a_{s^*y^*}+(1-u)(ua_{y^*}+(u+1)a_{s^*y^*})\\&=
-(u^2-u-1)a_{s^*y^*}-(u^2-u)a_{y^*}-(1-u^2)a_{s^*y^*}\endalign$$
which is obvious.

Assume next that $sw\ne ws^*>w$. We must show that $B(a_{sws^*})=T_s\i B(a_w)$ or that
$$\e_wT_{s^*w^*s}\i a_{s^*w^*s}=T_s\i\e_wT_{w^*}\i a_{w^*}$$
or that
$$T_s\i T_{w^*}\i T_{s^*}\i a_{s^*w^*s}=T_s\i T_{w^*}\i a_{w^*}$$
or that 
$$a_{s^*w^*s}=T_{s^*}a_{w^*}.$$
This follows from 0.1(iii) with $s,w$ replaced by $s^*,w^*$.

Finally assume that $sw\ne ws^*>w$. We set $y=sws^*\in\II_*$ so that $sy>y$. We must show that 
$B((u^2-1)a_{sys^*}+u^2a_y)=T_s\i B(a_{sys^*})$ or that
$$(u^{-2}-1)\e_yT_{s^*y^*s}\i a_{s^*y^*s}+u^{-2}\e_yT_{y^*}\i a_{y^*}=T_s\i \e_yT_{s^*y^*s}\i a_{s^*y^*s}$$
or that 
$$(u^{-2}-1)T_s\i T_{y^*}\i T_{s^*}\i a_{s^*y^*s}+u^{-2}T_{y^*}\i a_{y^*}
=T_s\i T_s\i T_{y^*}\i T_{s^*}\i a_{s^*y^*s}$$
or (using 0.1(iii) with $w,s$ replaced by $y^*,s^*$) that
$$(u^{-2}-1)T_s\i T_{y^*}\i a_{y^*}+u^{-2}T_{y^*}\i a_{y^*}=T_s\i T_s\i T_{y^*}\i a_{y^*}$$
or that
$$(u^{-2}-1)T_s\i+u^{-2}=T_s\i T_s\i$$
which is obvious.

This completes the proof of (a). Since the elements $T_s$ generate the algebra $\fH$, from (a) we deduce that
$B(hm)=\bar{h}B(m)$ for any $h\in\fH,m\in M$. This proves the existence part of 0.2(a).

For $n\in\ZZ,w\in\II_*$ we have
$$B(B(u^na_w))=\e_wB(u^{-n}T_{w^*}\i a_{w^*})=\e_w\e_{w^*}u^n T_{w^{*-1}}T_w\i a_w=u^na_w.$$
Thus $B^2=1$. The uniqueness part of 0.2(a) is proved as in \cite{\LV, 2.9}. This completes the proof of 0.2(a). 
Now 0.2(b) follows from the proof of 0.2(a).

\head 4. Proof of Theorem 0.4\endhead
\subhead 4.1\endsubhead
For $w\in\II_*$ we have
$$\ov{a'_w}=\sum_{y\in\II_*}\ov{r_{y,w}}a'_y$$
where $r_{y,w}\in\uca$ is zero for all but finitely many $y$. (This $r_{y,w}$ differs from that in 
\cite{\LV, 0.2(b)}.)

For $s\in S$ we set $T'_s=u\i T_s$. We rewrite the formulas 0.1(i)-(iv) as follows.

(i) $T'_sa'_w=a'_w+(v+v\i)a'_{sw}$ if $sw=ws^*>w$;

(ii) $T'_sa'_w=(u-1-u\i)a'_w+(v-v\i)a'_{sw}$ if $sw=ws^*<w$;

(iii) $T'_sa'_w=a'_{sws^*}$ if $sw\ne ws^*>w$;

(iv) $T'_sa'_w=(u-u\i)a'_w+a'_{sws^*}$ if $sw\ne ws^*<w$.

\subhead 4.2\endsubhead
Now assume that $y\in\II_*,sy>y$. From the equality $\ov{T'_sa'_y}=\ov{T'_s}(\ov{a'_y})$ (where  
$\ov{T'_s}=T'_s+u\i-u$) we see that
$$\sum_x\ov{r_{x,y}}a'_x+(v+v\i)\sum_x\ov{r_{x,sy}}a'_x\text{ (if $sy=ys^*$) or }\sum_x\ov{r_{x,sys^*}}a'_x
\text{ (if $sy\ne ys^*$)}$$
is equal to
$$\align&\sum_{x;sx=xs^*,sx>x}\ov{r_{x,y}}a'_x+\sum_{x;sx=xs^*,sx>x}\ov{r_{x,y}}(v+v\i)a'_{sx}\\&
+\sum_{x;sx=xs^*,sx<x}\ov{r_{x,y}}(u-1-u\i)a'_x+\sum_{x;sx=xs^*,sx<x}\ov{r_{x,y}}(v-v\i)a'_{sx}\\&
+\sum_{x;sx\ne xs^*,sx>x}\ov{r_{x,y}}a'_{sxs^*}
+\sum_{x;sx\ne xs^*,sx<x}\ov{r_{x,y}}(u-u\i)a'_x+\sum_{x;sx\ne xs^*,sx<x}\ov{r_{x,y}}a'_{sxs^*}\\&
+(u\i-u)\sum_x\ov{r_{x,y}}a'_x\\&
=
\sum_{x;sx=xs^*,sx>x}\ov{r_{x,y}}a'_x+\sum_{x;sx=xs^*,sx<x}\ov{r_{sx,y}}(v+v\i)a'_x\\&
+\sum_{x;sx=xs^*,sx<x}\ov{r_{x,y}}(u-1-u\i)a'_x+\sum_{x;sx=xs^*,sx>x}\ov{r_{sx,y}}(v-v\i)a'_x\\&
+\sum_{x;sx\ne xs^*,sx<x}\ov{r_{sxs^*,y}}a'_x
+\sum_{x;sx\ne xs^*,sx<x}\ov{r_{x,y}}(u-u\i)a'_x+\sum_{x;sx\ne xs^*,sx>x}\ov{r_{sxs^*,y}}a'_x\\&
+(u\i-u)\sum_x\ov{r_{x,y}}a'_x.\endalign$$
Hence when $sy=ys^*>y$ and $x\in\II_*$, we have
$$\align&(v+v\i)\ov{r_{x,sy}}=\ov{r_{sx,y}}(v-v\i)+(u\i-u)\ov{r_{x,y}}\text{ if }sx=xs^*>x,\\&
(v+v\i)\ov{r_{x,sy}}=-2\ov{r_{x,y}}+\ov{r_{sx,y}}(v+v\i)\text{ if }sx=xs^*<x,\\&
(v+v\i)\ov{r_{x,sy}}=\ov{r_{sxs^*,y}}+(u\i-1-u)\ov{r_{x,y}}\text{ if }sx\ne xs^*>x,\\&  
(v+v\i)\ov{r_{x,sy}}=-\ov{r_{x,y}}+\ov{r_{sxs^*,y}}\text{ if } sx\ne xs^*<x;\endalign$$
when $sy\ne ys^*>y$ and $x\in\II_*$, we have
$$\align&\ov{r_{x,sys^*}}=\ov{r_{sx,y}}(v-v\i)+(u\i+1-u)\ov{r_{x,y}}\text{ if } sx=xs^*>x,\\&
\ov{r_{x,sys^*}}=\ov{r_{sx,y}}(v+v\i)-\ov{r_{x,y}}\text{ if } sx=xs^*<x,\\&
\ov{r_{x,sys^*}}=\ov{r_{sxs^*,y}}+(u\i-u)\ov{r_{x,y}}\text{ if }sx\ne xs^*>x,\\&
\ov{r_{x,sys^*}}=\ov{r_{sxs^*,y}}\text{ if }sx\ne xs^*<x.\endalign$$
Applying $\,\bar{}$ we see that when $sy=ys^*>y$ and $x\in\II_*$, we have
$$\align&(v+v\i)r_{x,sy}=r_{sx,y}(v\i-v)+(u-u\i)r_{x,y}\text{ if }sx=xs^*>x,\\&
(v+v\i)r_{x,sy}=-2r_{x,y}+r_{sx,y}(v+v\i)\text{ if }sx=xs^*<x,\\&
(v+v\i)r_{x,sy}=r_{sxs^*,y}+(u-1-u\i)r_{x,y}\text{ if }sx\ne xs^*>x,\\&
(v+v\i)r_{x,sy}=-r_{x,y}+r_{sxs^*,y}\text{ if }sx\ne xs^*<x;\tag a\endalign$$
when $sy\ne ys^*>y$ and $x\in\II_*$, we have
$$\align&r_{x,sys^*}=r_{sx,y}(v\i-v)+(u+1-u\i)r_{x,y}\text{ if }sx=xs^*>x,\\&
r_{x,sys^*}=r_{sx,y}(v+v\i)-r_{x,y}\text{ if }sx=xs^*<x,\\&
r_{x,sys^*}=r_{sxs^*,y}+(u-u\i)r_{x,y}\text{ if } sx\ne xs^*>x,\\&
r_{x,sys^*}=r_{sxs^*,y}\text{ if }sx\ne xs^*<x.\tag b\endalign$$

\subhead 4.3\endsubhead
Setting $r'_{x,w}=v^{-l(w)+l(x)}r_{x,w}$, $r''_{x,w}=v^{-l(w)+l(x)}\ov{r_{x,w}}$ for $x,w\in\II_*$ we can 
rewrite the last formulas in 4.2 as follows.

When $x,y\in\II_*,sy=ys^*>y$ we have
$$\align&(v+v\i)vr'_{x,sy}=v\i r'_{sx,y}(v\i-v)+(u-u\i)r'_{x,y}\text{ if }sx=xs^*>x,\\&   
(v+v\i)vr'_{x,sy}=-2r'_{x,y}+r'_{sx,y}v(v+v\i)\text{ if } sx=xs^*<x,\\&              
(v+v\i)vr'_{x,sy}=v^{-2}r'_{sxs^*,y}+(u-1-u\i)r'_{x,y}\text{ if }sx\ne xs^*>x,\\&
(v+v\i)vr'_{x,sy}=-r'_{x,y}+v^2r'_{sxs^*,y}\text{ if } sx\ne xs^*<x.\endalign$$
When $x,y\in\II_*,sy\ne ys^*>y$, we have
$$\align&v^2r'_{x,sys^*}=r'_{sx,y}v\i(v\i-v)+(u+1-u\i)r'_{x,y}\text{ if }sx=xs^*>x,\\&
v^2r'_{x,sys^*}=r'_{sx,y}v(v+v\i)-r'_{x,y}\text{ if } sx=xs^*<x,\\&
v^2r'_{x,sys^*}=v^{-2}r'_{sxs^*,y}+(u-u\i)r'_{x,y}\text{ if }sx\ne xs^*>x,\\&
v^2r'_{x,sys^*}=v^2r'_{sxs^*,y}\text{ if }sx\ne xs^*<x.\endalign$$
When $x,y\in\II_*,sy=ys^*>y$ we have
$$\align&(v+v\i)vr''_{x,sy}=v\i r''_{sx,y}(v-v\i)+(u\i-u)r''_{x,y}\text{ if }sx=xs^*>x,\\&
(v+v\i)vr''_{x,sy}=-2r''_{x,y}+r''_{sx,y}v(v+v\i)\text{ if }sx=xs^*<x,\\&           
(v+v\i)vr''_{x,sy}=v^{-2}r''_{sxs^*,y}+(u\i-1-u)r''_{x,y}\text{ if } sx\ne xs^*>x,\\&
(v+v\i)vr''_{x,sy}=-r''_{x,y}+v^2r''_{sxs^*,y}\text{ if } sx\ne xs^*<x.\endalign$$
When $x,y\in\II_*,sy\ne ys^*>y$, we have
$$\align&v^2r''_{x,sys^*}=r''_{sx,y}v\i(v-v\i)+(u\i+1-u)r''_{x,y}\text{ if }sx=xs^*>x,\\&
v^2r''_{x,sys^*}=r''_{sx,y}v(v+v\i)-r''_{x,y}\text{ if }sx=xs^*<x,\\&
v^2r''_{x,sys^*}=v^{-2}r''_{sxs^*,y}+(u\i-u)r''_{x,y}\text{ if }sx\ne xs^*>x,\\&
v^2r''_{x,sys^*}=v^2r''_{sxs^*,y}\text{ if } sx\ne xs^*<x.\endalign$$

\proclaim{Proposition 4.4} Let $w\in\II_*$. 

(a) If $x\in\II_*,r_{x,w}\ne0$ then $x\le w$.

(b) If $x\in\II_*,x\le w$ we have $r'_{x,w}=\ZZ[v^{-2}]$, $r''_{x,w}=\ZZ[v^{-2}]$.
\endproclaim
We argue by induction on $l(w)$. If $w=1$ then $r_{x,w}=\d_{x,1}$ so that the result holds. Now assume that 
$l(w)\ge1$. We can find $s\in S$ such that $sw<w$. Let $y=s\bul w\in\II_*$ (see 0.6). We have $y<w$. In the setup
of (a) we have $r_{x,s\bul y}\ne0$. From the formulas in 4.3 we deduce the following.

If $sx=xs^*$ then $r'_{sx,y}\ne0$ or $r'_{x,y}\ne0$ hence (by the induction hypothesis) $sx\le y$ or $x\le y$; if
$x\le y$ then $x\le w$ while if $sx\le y$ we have $sx\le w$ hence by \cite{\HEC,2.5} we have $x\le w$.

If $sx\ne xs^*$ then $r'_{sxs^*,y}\ne0$ or $r'_{x,y}\ne0$ hence (by the induction hypothesis) $sxs^*\le y$ or 
$x\le y$; if $x\le y$ then $x\le w$ while if $sxs^*\le y$ we have $sxs^*\le w$ hence by \cite{\HEC, 2.5} we have 
$x\le w$.

We see that $x\le w$ and (a) is proved.

In the remainder of the proof we assume that $x\le w$. Assume that $sy=ys^*$. Using the formulas in 4.3 and the 
induction hypothesis we see that $v(v+v\i)r'_{x,w}\in v^2\ZZ[v^{-2}]$, $v(v+v\i)r''_{x,w}\in v^2\ZZ[v^{-2}]$;
hence $r'_{x,w}\in\ZZ[[v^{-2}]]$, $r''_{x,w}\in\ZZ[[v^{-2}]]$. Since $r'_{x,w}\in\ZZ[v,v\i]$,
$r''_{x,w}\in\ZZ[v,v\i]$, it follows that $r'_{x,w}\in\ZZ[v^{-2}]$, $r''_{x,w}\in\ZZ[v^{-2}]$. 

Assume now that  $sy\ne ys^*$. Using the formulas in 4.3 and the induction hypothesis we see that
$v^2r'_{x,w}\in v^2\ZZ[v^{-2}]$, $v^2r''_{x,w}\in v^2\ZZ[v^{-2}]$; hence $r'_{x,w}\in\ZZ[v^{-2}]$, 
$r''_{x,w}\in\ZZ[v^{-2}]$. This completes the proof.

\proclaim{Proposition 4.5} (a) There is a unique function $\ph:\II_*@>>>\NN$ such that $\ph(1)=0$ and for any 
$w\in\II_*$ and any $s\in S$ with $sw<w$ we have $\ph(w)=\ph(sw)+1$ (if $sw=ws^*$) and $\ph(w)=\ph(sws^*)$ (if 
$sw\ne ws^*$). For any $w\in\II_*$ we have $l(w)=\ph(w)\mod2$. Hence, setting $\k(w)=(-1)^{(l(w)+\ph(w))/2}$ for 
$w\in\II_*$ we have $\k(1)=1$ and $\k(w)=-\k(s\bul w)$ (see 0.6) for any $s\in S,w\in\II_*$ such that $sw<w$.

(b) If $x,w\in\II_*,x\le w$ then the constant term of $r'_{x,w}$ is $1$ and the constant term of $r''_{x,w}$ is 
$\k(x)\k(w)$ (see 4.4(b)).
\endproclaim
We prove (a). Assume first that $*$ is the identity map. For $w\in\II_*$ let $\ph(w)$ be the dimension of the 
$-1$ eigenspace of $w$ on the reflection representation of $W$. This function has the required properties. If $*$
is not the identity map, the proof is similar: for $w\in\II_*$, $\ph(w)$ is the dimension of the $-1$ eigenspace 
of $wT$ minus the dimension of the $-1$ eigenspace of $T$ where $T$ is an automorphism of the reflection 
representation of $W$ induced by $*$. 

We prove (b). Let $n'_{x,w}$ (resp. $n''_{x,w}$) be the constant term of $r'_{x,w}$ (resp. $r''_{x,w}$). We shall 
prove for any $w\in\II_*$ the following statement:

(c) {\it If $x\in\II_*,x\le w$ then $n'_{x,w}=1$ and $n''_{x,w}=n''_{1,x}n''_{1,w}\in\{1,-1\}$.}
\nl
We argue by induction on $l(w)$. If $w=1$ we have $r'_{w,w}=r''_{w,w}=1$ and (c) is obvious. We assume that 
$w\in\II_*,w\ne 1$. We can find $s\in S$ such that $sw<w$. We set $y=s\bul w$. Taking the coefficients of $v^2$ in
the formulas in 4.3 and using 4.4(b) we see that the following holds for any $x\in\II_*$ such that $x\le w$:
$$n'_{x,w}=n'_{x,y}, n''_{x,w}=-n''_{x,y}\text{ if }sx>x,$$
(by \cite{\HEC, 2.5(b)}, we must have $x\le y$) and
$$n'_{x,w}=n'_{s\bul x,y}, n''_{x,w}=n''_{s\bul x,y}\text{ if }sx<x$$
(by \cite{\HEC, 2.5(b)}, we must have $s\bul x\le y$).
\nl
Using the induction hypothesis we see that $n'_{x,w}=1$ and
$$n''_{x,w}=-n''_{1,x}n''_{1,y}\text{ if }sx>x,$$          
$$n''_{x,w}=n''_{1,s\bul x}n''_{1,y}\text{ if }sx<x.$$
Also, taking $x=1$ we see that 
$$n''_{1,w}=-n''_{1,y}.\tag d$$
Returning to a general $x$ we deduce
$$n''_{x,w}=n''_{1,x}n''_{1,w}\text{ if }sx>x,$$          
$$n''_{x,w}=-n''_{1,s\bul x}n''_{1,w}\text{ if }sx<x.$$
Applying (d) with $w$ replaced by $x$ we see that $n''_{1,x}=-n''_{1,s\bul x}$ if $sx<x$. This shows by
induction on $l(x)$ that $n''_{1,x}=\k(x)$ for any $x\in\II_*$. Thus we have 
$n''_{x,w}=n''_{1,x}n''_{1,w}=\k(x)\k(w)$ for any $x\le w$. This completes the inductive proof of (c) and that 
of (b). The proposition is proved.

\subhead 4.6\endsubhead
We show:

(a) {\it For any $x,z\in\II_*$ such that $x\le z$ we have 
$\sum_{y\in\II_*;x\le y\le z}\ov{r_{x,y}}r_{y,z}=\d_{x,z}$.}
\nl
Using the fact that $\,\bar{}\,:uM@>>>\uM$ is an involution we have
$$a'_z=\ov{\ov{a'_z}}=\ov{\sum_{y\in\II_*}\ov{r_{y,z}a'_y}}
=\sum_{y\in\II_*}r_{y,z}\ov{a'_y}=\sum_{y\in\II_*}\sum_{x\in\II_*}r_{y,z}\ov{r_{x,y}}a'_x.$$
We now compare the coefficients of $a'_x$ on both sides and use 4.4(a); (a) follows.

The following result provides the M\"obius function for the partially ordered set $(\II_*,\le)$.

\proclaim{Proposition 4.7} Let $x,z\in\II_*,x\le z$. Then $\sum_{y\in\II_*;x\le y\le z}\k(x)\k(y)=\d_{x,z}$.
\endproclaim
We can assume that $x<z$. By 4.4(b), 4.5(b) for any $y\in\II_*$ such that $x\le y\le z$ we have 
$$\ov{r_{x,y}}r_{y,z}=v^{l(y)-l(x)}v^{l(z)-l(x)}r''_{x,y}r'_{y,z}\in v^{l(z)-l(x)}(\k(x)\k(y)+v^{-2}\ZZ[v^{-2}]).
$$
Hence the identity 4.6(a) implies that
$$\sum_{y\in\II_*;x\le y\le z}v^{l(z)-l(x)}\k(x)\k(y)+\text{ strictly lower powers of  $v$ is $0$}.$$
In particular, $\sum_{y\in\II_*;x\le y\le z}\k(x)\k(y)=0$. The proposition is proved.

\subhead 4.8\endsubhead
For any $w\in\II_*$ we have
$$r_{w,w}=1.\tag a$$
Indeed by 4.4(b) we have $r_{w,w}\in\ZZ[v^{-2}]$, $\ov{r_{w,w}}\in\ZZ[v^{-2}]$ hence $r_{w,w}$ is a constant. By 
4.5(b) this constant is $1$.

\subhead 4.9\endsubhead
Let $w\in\II_*$. We will construct for any $x\in\II_*$ such that $x\le w$ an element $u_x\in\uca_{\le 0}$ such 
that

(a) $u_x=1$,

(b) $u_x\in\uca_{<0}$, $\ov{u_x}-u_x=\sum_{y\in\II_*;x<y\le w}r_{x,y}u_y$ for any $x<w$.
\nl
The argument is almost a copy of one in \cite{\HEC, 5.2}. We argue by induction on $l(w)-l(x)$. If $l(w)-l(x)=0$ 
then $x=w$ and we set $u_x=1$. Assume now that $l(w)-l(x)>0$ and that $u_z$ is already defined whenever
$z\le w$, $l(w)-l(z)<l(w)-l(x)$ so that (a) holds and (b) holds if $x$ is replaced by any such $z$.
Then the right hand side of the equality in (b) is defined. We denote it by $\a_x\in\uca$. We have
$$\align&\a_x+\bar\a_x=\sum_{y\in\II_*;x<y\le w}r_{x,y}u_y+\sum_{y\in\II_*;x<y\le w}\ov{r_{x,y}}\bar u_y\\&
=\sum_{y\in\II_*;x<y\le w}r_{x,y}u_y
+\sum_{y\in\II_*;x<y\le w}\ov{r_{x,y}}(u_y+\sum_{z\in\II_*;y<z\le w}r_{y,z}u_z)\\&
=\sum_{y\in\II_*;x<y\le w}r_{x,y}u_y+\sum_{z\in\II_*;x<z\le w}\ov{r_{x,z}}u_z
+\sum_{z\in\II_*;x<z\le w}\sum_{y\in\II_*;x<y<z}\ov{r_{x,y}}r_{y,z}u_z\\&
=\sum_{z\in\II_*;x<z\le w}\sum_{y\in\II_*;x\le y<z}\ov{r_{x,y}}r_{y,z}u_z=\sum_{z\in\II_*;x<z\le w}\d_{x,z}u_z=0.
\endalign$$
(We have used 4.6(a), 4.8(a).) Since $\a_x+\bar\a_x=0$ we have $\a_x=\sum_{n\in\ZZ}\g_nv^n$ (finite sum) where
$\g_n\in\ZZ$ satisfy $\g_n+\g_{-n}=0$ for all $n$ and in particular $\g_0=0$. Then
$u_x=-\sum_{n<0}\g_nv^n\in\uca_{<0}$ satisfies $\bar u_x-u_x=\a_x$. This completes the inductive construction of 
the elements $u_x$.

We set $A_w=\sum_{y\in\II_*;y\le w}u_ya'_y\in\uM_{\le 0}$. We have
$$\align&\ov{A_w}=\sum_{y\in\II_*;y\le w}\bar u_y\ov{a'_y}=\sum_{y\in\II_*;y\le w}\bar u_y
\sum_{x\in\II_*;x\le y}\ov{r_{x,y}}a'_x\\&=\sum_{x\in\II_*;x\le w}(\sum_{y\in\II_*;x\le y\le w}\ov{r_{x,y}}
\bar u_y)a'_x=\sum_{x\in\II_*;x\le w}u_xa'_x=A_w.\endalign$$
We will also write $u_y=\p_{y,w}\in\uca_{\le0}$ so that 
$$A_w=\sum_{y\in\II_*;y\le w}\p_{y,w}a'_y.$$
Note that $\p_{w,w}=1$, $\p_{y,w}\in\uca_{<0}$ if $y<w$ and
$$\ov{\p_{y,w}}=\sum_{z\in\II_*;y\le z\le w}r_{y,z}\p_{z,w}.$$
We show that for any $x\in\II_*$ such that $x\le w$ we have:

(c) {\it $v^{l(w)-l(x)}\p_{x,w}\in\ZZ[v]$ and has constant term $1$.}
\nl
We argue by induction on $l(w)-l(x)$. If $l(w)-l(x)=0$ then $x=w$, $\p_{x,w}=1$ and the result is obvious. Assume
now that $l(w)-l(x)>0$. Using 4.4(b) and 4.5(b) and the induction hypothesis we see that
$$\sum_{y\in\II_*;x<y\le w}r_{x,y}\p_{y,w}=\sum_{y\in\II_*;x<y\le w}v^{-l(y)+l(x)}\ov{r''_{x,y}}\p_{y,w}$$
is equal to
$$\sum_{y\in\II_*;x<y\le w}v^{-l(y)+l(x)}\k(x)\k(y)v^{-l(w)+l(y)}=v^{-l(w)+l(x)}\sum_{y\in\II_*;x<y\le w}
\k(x)\k(y)$$
plus strictly higher powers of $v$. Using 4.7, this is $-v^{-l(w)+l(x)}$ plus strictly higher powers of $v$. Thus,
$$\ov{\p_{x,w}}-\p_{x,w}=-v^{-l(w)+l(x)}+\text{ plus strictly higher powers of }v.$$
Since $\ov{\p_{x,w}}\in v\ZZ[v]$, it is in particular a $\ZZ$-linear combination of powers of $v$
strictly higher than $-l(w)+l(x)$. Hence
$$-\p_{x,w}=-v^{-l(w)+l(x)}+\text{ plus strictly higher powers of }v.$$
This proves (c).

We now show that for any $x\in\II_*$ such that $x\le w$ we have:
$$v^{l(w)-l(x)}\p_{x,w}\in\ZZ[u,u\i].\tag d$$
We argue by induction on $l(w)-l(x)$. If $l(w)-l(x)=0$ then $x=w$, $\p_{x,w}=1$ and the result is obvious. Assume
now that $l(w)-l(x)>0$. Using 4.4(b) and the induction hypothesis we see that 
$$\sum_{y\in\II_*;x<y\le w}r_{x,y}\p_{y,w}=\sum_{y\in\II_*;x<y\le w}v^{-l(y)+l(x)}\ov{r''_{x,y}}\p_{y,w}$$
belongs to 
$$\sum_{y\in\II_*;x<y\le w}v^{-l(y)+l(x)}v^{-l(w)+l(y)}\ZZ[v^2,v^{-2}]$$ 
hence to $v^{-l(w)+l(x)}\ZZ[v^2,v^{-2}]$. Thus,
$$\ov{\p_{x,w}}-\p_{x,w}\in v^{-l(w)+l(x)}\ZZ[v^2,v^{-2}].$$
It follows that both $\ov{\p_{x,w}}$ and $\p_{x,w}$ belong to $v^{-l(w)+l(x)}\ZZ[v^2,v^{-2}]$. This proves (d).

Combining (c), (d) we see that for any $x\in\II_*$ such that $x\le w$ we have:

(e) {\it $v^{l(w)-l(x)}\p_{x,w}=P^\pm_{x,w}$ where $P^\pm_{x,w}\in\ZZ[u]$ has constant term $1$.}
\nl
We have 
$$A_w=v^{-l(w)}\sum_{y\in\II_*;y\le w}P^\pm_{y,w}a_y.$$
Also, $P^\pm_{w,w}=1$ and for any $y\in\II_*$, $y<w$, we have $\deg P^\pm_{y,w}\le(l(w)-l(y)-1)/2$ (since
$\p_{y,w}\in\uca_{<0}$). Thus the existence statement in 0.4(a) is established. To prove the uniqueness statement
in 0.4(a) it is enough to prove the following statement:

(f) {\it Let $m,m'\in\uM$ be such that $\bar m=\bar m'$, $m-m'\in\uM_{>0}$. Then $m=m'$.}
\nl
The proof is entirely similar to that in \cite{\LV, 3.2} (or that of \cite{\HEC, 5.2(e)}). The proof of 0.4(b) is
immediate. This completes the proof of Theorem 0.4.

The following result is a restatement of (e).

\proclaim{Proposition 4.10} Let $y,w\in\II_*$ be such that $y\le w$. The constant term of $P_{y,w}^\pm\in\ZZ[u]$ 
is equal to $1$.
\endproclaim 

\head 5. The submodule $\uM^K$ of $\uM$\endhead
\subhead 5.1\endsubhead
Let $K$ be a subset of $S$ which generates a finite subgroup $W_K$ of $W$ and let $K^*$ be the image of $K$ under
$*$. For any $(W_K,W_{K^*})$-double coset $\Om$ in $W$ we denote by $d_\Om$ (resp. $b_\Om$) the unique element 
of maximal (resp. minimal) length of $\Om$. Now $w\m w^{*-1}$ maps any $(W_K,W_{K^*})$-double coset in $W$ to a 
$(W_K,W_{K^*})$-double coset in $W$; let $\II_*^K$ be the set of $(W_K,W_{K^*})$-double cosets $\Om$ in $W$ such 
that $\Om$ is stable under this map, or equivalently, such that $d_\Om\in\II_*$, or such that $b_\Om\in\II_*$.
We set 
$$\PP_K=\sum_{x\in W_K}u^{l(x)}\in\NN[u].$$
If in addition $K$ is $*$-stable we set
$$\PP_{H,*}=\sum_{x\in W_K,x^*=x}u^{l(x)}\in\NN[u].$$

\proclaim{Lemma 5.2} Let $\Om\in\II^K_*$. Let $x\in\II_*\cap\Om$ and let $b=b_\Om$. Then there exists a sequence
$x=x_0,x_1,\do,x_n=b$ in $\II_*\cap\Om$ and a sequence $s_1,s_2,\do,s_n$ in $S$ such that for any $i\in[1,n]$ we
have $x_i=s_i\bul x_{i-1}$.
\endproclaim
We argue by induction on $l(x)$ (which is $\ge l(b)$). If $l(x)=l(b)$ then $x=b$ and the result is obvious (with 
$n=0$). Now assume that $l(x)>l(b)$. Let $H=K\cap(bK'b\i)$. By 1.2(a) we have $x=cbzc^{*-1}$ where 
$c\in W_K$, $z\in W_{H^*}$ satisfies $bz=z^*b$ and $l(x)=l(c)+l(b)+l(z)+l(c)$. If $c\ne1$ we write $c=sc'$, 
$s\in K,c'\in W_K$, $c'<c$ and we set $x_1=c'bzc'{}^{*-1}$. We have $x_1=sxs^*\in\Om$, $l(x_1)<l(x)$. Using the 
induction hypothesis for $x_1$ we see that the desired result holds for $x$. Thus we can assume that $c=1$ so 
that $x=bz$. Let $\t:W_{H^*}@>>>W_{H^*}$ be the automorphism $y\m b\i y^*b$; note that $\t(H^*)=H^*$ and 
$\t^2=1$. We have $z\in\II_\t$ where $\II_\t:=\{y\in W_{H^*};\t(y)\i=y\}$. 

Since $l(bz)>l(b)$ we have $z\ne1$. We can find $s\in H^*$ such that $sz<z$. 

If $sz=z\t(s)$ then $sz\in\II_\t$, $bsz\in\Om$, $l(bsz)<l(bz)$. Using the induction hypothesis for $bsz$ instead
of $x$ we see that the desired result holds for $x=bz$. (We have $bsz=tbz=bzt^*$ where $t=(\t(s))^*\in H$.)

If $sz\ne z\t(s)$ then $sz\t(s)\in\II_t$, $bsz\t(s)\in\Om$, $l(bsz\t(s))<l(bz)$. Using the induction hypothesis 
for $bsz\t(s)$ instead of $x$ we see that the desired result holds for $x=bz$. (We have $bsz\t(s)=tbzt^*$ where 
$t=(\t(s))^*\in H$.) The lemma is proved.

\subhead 5.3\endsubhead
For any $\Om\in\II_*^K$ we set 
$$a_\Om=\sum_{w\in\II_*\cap\Om}a_w\in\uM.$$
Let $\uM^K$ be the $\uca$-submodule of $\uM$ spanned by the elements $a_\Om (\Om\in\II^K_*)$. In other words, 
$\uM^K$ consists of all $m=\sum_{w\in\II_*}m_wa_w\in\uM$ such that the function $\II_*@>>>\uca$ given by 
$w\m m_w$ is constant on $\II_*\cap\Om$ for any $\Om\in\II_*$.

\proclaim{Lemma 5.4} (a) We have $\uM^K=\cap_{s\in K}\uM^{\{s\}}$.

(b) The $\uca$-submodule $\uM^K$ is stable under $\,\bar{}:\uM@>>>\uM$.

(c) Let $\SS=\sum_{x\in W_K}T_x\in\ufH$ and let $m\in\uM$. We have $\SS m\in\uM^K$.
\endproclaim
We prove (a). The fact that $\uM^K\sub\uM^{\{s\}}$ (for $s\in K$) follows from the fact that any 
$(W_K,W_{K^*})$-double coset in $W$ is a union of $(W_{\{s\}},W_{\{s^*\}})$-double cosets in $W$. Thus we have 
$\uM^K\sub\cap_{s\in K}\uM^{\{s\}}$. Conversely let $m\in\cap_{s\in K}\uM^{\{s\}}$. We have
$m=\sum_{w\in\II_*}m_wa_w\in\uM$ where $m_w\in\uca$ is zero for all but finitely many $w$ and we have 
$m_w=m_{s\bul w}$ if $w\in\II_*,s\in K$. Using 5.2 we see 
that $m_x=m_{b_\Om}=m_{x'}$ whenever $x,x'\in\II_*$ are in the same $(W_K,W_{K^*})$-double coset $\Om$ in $W$.
Thus, $m\in\uM^K$. This proves (a).

We prove (b). Using (a), we can assume that $K=\{s\}$ with $s\in S$. By 1.3, if $\Om\in\II_*^{\{s\}}$, then we 
have
$\Om=\{w,s\bul w\}$ for some $w\in\II_*$ such that $sw>w$. Hence it is enough to show that for such $w$ we have
$\ov{a_w+a_{s\bul w}}\in\uM^{\{s\}}$. We have $\ov{a_w+a_{s\bul w}}=\sum_{x\in\II_*}m_xa_x$ with $m_x\in\uca$ and
we must show that $m_x=m_{s\bul x}$ for any $x\in\II_*$. If we can show that 
$f\ov{a_w+a_{s\bul w}}\in\uM^{\{s\}}$ for some $f\in\uca-\{0\}$ then it would follow that for any $x\in\II_*$ we 
have $fm_x=fm_{s\bul x}$ hence $m_x=m_{s\bul x}$ as desired. Thus it is enough to show that  

(d) $(u\i+1)\ov{a_w+a_{sw}}\in\uM^{\{s\}}$ if $w\in\II_*$ is such that $sw=ws^*>w$,

(e) $\ov{a_w+a_{sws^*}}\in\uM^{\{s\}}$ if $w\in\II_*$ is such that $sw\ne ws^*>w$.
\nl
In the setup of (d) we have 
$$\align&(u\i+1)\ov{a_w+a_{sw}}=\ov{(u+1)(a_w+a_{sw})}
=\ov{(T_s+1)a_w}=\ov{T_s+1}(\ov{a_w})\\&=u^{-2}(T_s+1)\ov{a_w}\endalign$$
(see 0.1(i)); in the setup of (e) we have  
$$\ov{a_w+a_{sws^*}}=\ov{(T_s+1)a_w}=\ov{T_s+1}(\ov{a_w})=u^{-2}(T_s+1)(\ov{a_w})$$
(see 0.1(iii)). Thus it is enough show that $(T_s+1)(\ov{a_w})\in\uM^{\{s\}}$ for any $w\in\II_*$. Since 
$\ov{a_w}$ is an $\uca$-linear combination of elements $a_x,x\in\II_*$ it is enough to show that
$(T_s+1)a_x\in\uM^{\{s\}}$. This follows immediately from 0.1(i)-(iv). 

We prove (c). Let $m'=\SS m=\sum_{w\in\II_*}m'_wa_w$, $m'_w\in\uca$. For any $s\in K$ we have 
$\SS=(T_s+1)h$ for some $h\in\ufH$ hence $m'\in(T_s+1)\uM$. This implies by the formulas 0.1(i)-(iv) that
$m'_w=w'_{s\bul w}$ for any $w\in\II_*$; in other words we have $m'\in\uM^{\{s\}}$. Since this holds for any 
$s\in K$ we see, using (a), that $m'\in\uM^K$. 
The lemma is proved.

\subhead 5.5\endsubhead
For $\Om,\Om'\in\II^K_*$ we write $\Om\le\Om'$ when $d_\Om\le d_{\Om'}$. This is a partial order on $\II^K_*$.
For any $\Om\in\II^K_*$ we set 
$$a'_\Om=v^{-l(d_\Om)}a_\Om=\sum_{x\in\Om\cap\II^K_*}v^{l(x)-l(d_\Om)}a'_x.$$
Clearly, $\{a'_{\Om'};\Om'\in\II^K_*\}$ is an $\uca$-basis of $\uM^K$. Hence from 5.4(b) we see that 
$$\ov{a'_{\Om}}=\sum_{\Om'\in\II^K_*}\ov{r_{\Om',\Om}}a'_{\Om'}$$
where $r_{\Om',\Om}\in\uca$ is zero for all but finitely many $\Om'$. On the other hand we have
$$\ov{a'_{\Om}}=\sum_{x\in\Om\cap\II_*,y\in\II_*;y\le x}v^{-l(x)+l(d_\Om)}\ov{r_{y,x}}a'_y\tag a$$
hence
$$r_{\Om',\Om}=\sum_{x\in\Om\cap\II_*;d_{\Om'}\le x}v^{l(x)-l(d_\Om)}r_{d_{\Om'},x}$$ 
It follows that
$$r_{\Om,\Om}=1\tag b$$
(we use that $r_{d_\Om,d_\Om}=1$) and
$$r_{\Om',\Om}\ne0\implies\Om'\le\Om.\tag c$$
Indeed, if for some $x\in\Om\cap\II_*$ we have $d_{\Om'}\le x$, then $d_{\Om'}\le d_\Om$. We have
$$a'_\Om=\ov{\ov{a'_\Om}}=\ov{\sum_{\Om'\in\II^K_*}\ov{r_{\Om',\Om}}a'_{\Om'}}=
\sum_{\Om'\in\II^K_*}r_{\Om',\Om}\sum_{\Om''\in\II^K_*}\ov{r_{\Om'',\Om'}}a'_{\Om''}.$$
Hence 
$$\sum_{\Om'\in\II^K_*}\ov{r_{\Om'',\Om'}}r_{\Om',\Om}=\d_{\Om,\Om''}\tag d$$
for any $\Om,\Om''$ in $\II^K_*$.

Note that
$$a'_\Om=a'_{d_\Om}\mod\uM_{<0}.\tag e$$
Indeed, if $x\in\Om\cap\II^K_*$, $x\ne d_\Om$ then $l(x)-l(d_\Om)<0$.

\subhead 5.6\endsubhead
Let $\Om\in\II_*^K$. We will construct for any $\Om'\in\II_*^K$ such that $\Om'\le\Om$ an element 
$u_{\Om'}\in\uca_{\le 0}$ such that

(a) $u_\Om=1$,

(b) $u_{\Om'}\in\uca_{<0}$, 
$\ov{u_{\Om'}}-u_{\Om'}=\sum_{\Om''\in\II_*^K;\Om'<\Om''\le\Om}r_{\Om',\Om''}u_{\Om''}$ for any $\Om'<\Om$.
\nl
The proof follows closely that in 4.9. We argue by induction on $l(d_\Om)-l(d_{\Om'})$. If 
$l(d_\Om)-l(d_{\Om'})=0$ then $\Om=\Om'$ and we set $u_{\Om'}=1$. 
Assume now that $l(d_\Om)-l(d_{\Om'})>0$ and that $u_{\Om_1}$ is already defined whenever
$\Om_1\le\Om$, $l(d_\Om)-l(d_{\Om_1})<l(d_\Om)-l(d_{\Om'})$ so that (a) holds and (b) holds if $\Om'$ is 
replaced by any such $\Om_1$.
Then the right hand side of the equality in (b) is defined. We denote it by $\a_{\Om'}\in\uca$. We have
$\a_{\Om'}+\ov{\a_{\Om'}}=0$ by a computation like that in 4.9, but using 5.5(b),(c),(d). From this we see that 
$\a_{\Om'}=\sum_{n\in\ZZ}\g_nv^n$ (finite sum) where $\g_n\in\ZZ$ satisfy $\g_n+\g_{-n}=0$ for all $n$ and in 
particular $\g_0=0$. Then $u_{\Om'}=-\sum_{n<0}\g_nv^n\in\uca_{<0}$ satisfies $\ov{u_{\Om'}}-u_{\Om'}=\a_{\Om'}$.
This completes the inductive construction of the elements $u_{\Om'}$.

We set $A_\Om=\sum_{\Om'\in\II^K_*;\Om'\le\Om}u_{\Om'}a'_{\Om'}\in\uM_{\le 0}\cap\uM^K$. We have 
$$\ov{A_\Om}=A_\Om.\tag c$$
(This follows from (b) as in the proof of the analogous equality 
$\ov{A_w}=A_w$ in 4.9.) We will also write $u_{\Om'}=\p_{\Om',\Om}\in\uca_{\le0}$ so that 
$$A_\Om=\sum_{\Om'\in\II_*^K;\Om'\le\Om}\p_{\Om',\Om}a'_{\Om'}.$$
We show 
$$A_\Om-A_{d_\Om}\in\uM_{<0}.\tag d$$
Using 5.5(a) and $\p_{\Om',\Om}\in\uca_{<0}$ (for $\Om'<\Om)$ we see that $A_\Om=a'_{d_\Om}\mod\uM_{<0}$; it 
remains to use that $A_{d_\Om}=a'_{d_\Om}\mod\uM_{<0}$. 

Applying 4.9(f) to $m=A_\Om$, $m'=A_{d_\Om}$ (we use (c),(d)) we deduce:
$$A_\Om=A_{d_\Om}.\tag e$$
In particular,

(f) {\it For any $\Om\in\II^K_*$ we have $A_{d_\Om}\in\uM^K$.}

\subhead 5.7\endsubhead 
We define an $\ca$-linear map $\z:M@>>>\QQ(u)$ by $\z(a_w)=u^{l(w)}(\fra{u-1}{u+1})^{\ph(w)}$ (see 4.5(a)) 
for $w\in\II_*$. We show:

(a) {\it For any $x\in W,m\in M$ we have $\z(T_xm)=u^{2l(x)}\z(m)$.}
\nl
We can assume that $x=s,m=a_w$ where $s\in S,w\in\II_*$. Then we are in 
one of the four cases (i)-(iv) in 0.1. We set $n=l(w)$, $d=\ph(w)$, $\lambda=\fra{u-1}{u+1}$. The identities to 
be checked in the cases 0.1(i)-(iv) are:
$$u^2u^n\l^d=uu^n\l^d+(u+1)u^{n+1}\l^{d+1},$$
$$u^2u^n\l^d=(u^2-u-1)u^n\l^d+(u^2-u)u^{n-1}\l^{d-1},$$
$$u^2u^n\l^d=u^{n+2}\l^d,$$
$$u^2u^n\l^d=(u^2-1)u^n\l^d+u^2u^{n-2}\l^d,$$ 
respectively. These are easily verified.

\subhead 5.8\endsubhead 
Assuming that $K^*=K$, we set 
$$\car_{K,*}=\sum_{y\in W_K;y^*=y\i}u^{l(y)}(\fra{u-1}{u+1})^{\ph(y)}\in\QQ(u).$$
Let $\Om\in\II^K_*$. Define $b,H,\t$ as in 5.2. Let 
$$W_K^H=\{c\in W_K;l(w)\le l(wr)\text{ for any }r\in W_H\}.$$
Using 1.2(a) we have
$\sum_{w\in\Om\cap\II_*}\z(a_w)=\sum_{c\in W_K^H}u^{2l(c)}\z(a_b)\car_{H^*,\t}(u)$ hence
$$\sum_{w\in\Om\cap\II_*}\z(a_w)=\PP_K(u^2)\PP_H(u^2)\i\z(a_b)\car_{H^*,\t}(u).\tag a$$
We have the following result.

\proclaim{Proposition 5.9} Assume that $W$ is finite. We have
$$\car_{S,*}(u)=\PP_S(u^2)\PP_{S,*}(u)\i.\tag a$$
\endproclaim
We can assume that $W$ is irreducible. We prove (a) by induction on $|S|$. If $|S|\le2$, (a) is easily checked. 
Now assume that $|S|\ge3$. Taking sum over all $\Om\in\II_*^K$ in 5.7(a) we obtain
$$\car_{S,*}(u)=\PP_K(u^2)\sum_{\Om\in\II^K_*}\PP_H(u^2)\i\z(a_b)\car_{H^*,\t}(u)$$
where $b,H,\t$ depend on $\Om$ as in 5.2. Using the induction hypothesis we obtain
$$\car_{S,*}(u)=\PP_K(u^2)\sum_{\Om\in\II^K_*}\z(a_b)\PP_{H^*,\t}(u)\i.$$
We now choose $K\sub S$ so that $W_K$ is of type 

$A_{n-1},B_{n-1},D_{n-1},A_1,B_3,A_5,D_7,E_7,I_2(5),H_3$ 
\nl
where $W$ is of type 

$A_n,B_n,D_n,G_2,F_4,E_6,E_7,E_8,H_3,H_4$
\nl
respectively. Then there are few $(W_K,W_{K^*})$ double cosets and the sum above can be computed in each case and
gives the desired result. 
(In the case where $W$ is a Weyl group, there is an alternative, uniform, proof of (a)
using flag manifolds over a finite field.) 

\subhead 5.10\endsubhead 
We return to the general case. Let $\Om\in\II^K_*$ and let $b,H,\t$ be as in 5.2. By 5.4(c) we have 
$\SS a_b\in\uM^K$. From 0.1(i)-(iv) we see that $\SS a_b=\sum_{y\in\Om\cap\II_*}f_ya_y$ where $f_y\in\ZZ[u]$ for 
all $y$. Hence we must have $\SS a_b=fa_\Om$ for some $f\in\ZZ[u]$. Appplying $\z$ to the last equality and using
5.7(a) we obtain
$\PP_K(u^2)\z(a_b)=f\sum_{y\in\Om\cap\II_*}\z(a_y)$. From 5.8(a), 5.9(a) we have
$$\sum_{y\in\Om\cap\II_*}\z(a_y)=\PP_K(u^2)\z(a_b)\PP_{H^*,\t}(u)\i$$
where $b,H,\t$ depend on $\Om$ as in 5.8. Thus $f=\PP_{H^*,\t}(u).$ We see that
$$\SS a_b=\PP_{H^*,\t}(u)a_\Om.\tag a$$

\subhead 5.11\endsubhead 
In this subsection we assume that $K^*=K$. Then $\Om:=W_K\in\II^K_*$. We have the following result.
$$A_\Om=v^{-l(w_K)}a_\Om.\tag a$$
By 5.6(f) we have $A_\Om=fa_\Om$ for some $f\in\uca$. Taking the coefficient of $a_{w_K}$ in both sides we get
$f=v^{-l(w_K)}$ proving (a).

Here is another proof of (a). It is enough to prove that $v^{-l(w_K)}a_\Om$ is fixed by 
$\,\bar{}$. By 5.10(a) we have $u^{-l(w_K)}\SS a_1=u^{-l(w_K)}\PP_{K,*}(u)a_\Om$. The left hand side of this 
equality is fixed by $\bar{}$ since $a_1$ and $u^{-l(w_K)}\SS$ are fixed by $\,\bar{}$.
Hence $v^{-2l(w_K)}\PP_{K,*}(u)a_\Om$ is fixed by $\,\bar{}$. 
Since $v^{-l(w_K)}\PP_{K,*}(u)$ is fixed by $\,\bar{}$ and is nonzero, it follows that
$v^{-l(w_K)}a_\Om$ is fixed by $\,\bar{}$, as desired.

\head 6. The action of $u\i(T_s+1)$ in the basis $(A_w)$ \endhead
\subhead 6.1\endsubhead
In this section we fix $s\in S$. 

Let $y,w\in\II_*$. When $y\le w$ we have as in 4.9, $\p_{y,w}=v^{-l(w)+l(y)}P^\pm_{y,w}$ so that 
$\p_{y,w}\in\uca_{<0}$ if $y<w$ and $\p_{w,w}=1$; when $y\not\le w$ we set $\p_{y,w}=0$. In any case we set as in
\cite{\LV, 4.1}:

(a) $\p_{y,w}=\d_{y,w}+\mu'_{y,w}v\i+\mu''_{y,w}v^{-2}\mod v^{-3}\ZZ[v\i]$
\nl
where $\mu'_{y,w}\in\ZZ,\mu''_{y,w}\in\ZZ$. Note that

(b) $\mu'_{y,w}\ne0\implies y<w, \e_y=-\e_w$,

(c) $\mu''_{y,w}\ne0 \implies y<w, \e_y=\e_w$.

\subhead 6.2\endsubhead
As in \cite{\LV, 4.3}, for any $y,w\in\II_*$ such that $sy<y<sw>w$ we define $\cm^s_{y,w}\in\uca$ by:
$$\cm^s_{y,w}=\mu''_{y,w}-\sum_{x\in\II_*;y<x<w,sx<x}\mu'_{y,x}\mu'_{x,w}-\d_{sw,ws^*}\mu'_{y,sw}
+\mu'_{sy,w}\d_{sy,ys^*}$$
if $\e_y=\e_w$,
$$\cm^s_{y,w}=\mu'_{y,w}(v+v\i)$$ 
if $\e_y=-\e_w$. 

The following result was proved in \cite{\LV, 4.4} assuming that $W$ is a Weyl group or affine Weyl group.
(We set $c_s=u\i(T_s+1)\in\ufH$.)

\proclaim{Theorem 6.3} Let $w\in\II_*$. 

(a) If $sw=ws^*>w$ then $c_sA_w=(v+v\i)A_{sw}+\sum_{z\in\II_*;sz<z<sw}\cm^s_{z,w}A_z$.

(b) If $sw\ne ws^*>w$ then $c_sA_w=A_{sws^*}+\sum_{z\in\II_*;sz<z<sws^*}\cm^s_{z,w}A_z$.

(c) If $sw<w$ then $c_sA_w=(u+u\i)A_w$.
\endproclaim
(In the case considered in \cite{\LV, 4.4} the last sum in the formula which corresponds to (b) involves 
$sz<z<sw$ instead of $sz<z<sws^*$; but as shown in {\it loc.cit.} the two conditions are equivalent.)

We prove (c). We have $sw<w$. By 5.6(f) we have $A_w\in\uM^{\{s\}}$. Hence it is enough to show that
$c_sm=(u+u\i)m$ where $m$ runs through a set of generators of the $\uca$-module $\uM^{\{s\}}$. Thus it is enough 
to show that $c_s(a_x+a_{s\bul x})=(u+u\i)(a_x+a_{s\bul x})$ for any $x\in\II_*$. This follows immediately from
0.1(i)-(iv).

Now the proof of (a),(b) (assuming (c)) is exactly as in \cite{\LV, 4.4}. (Note that in \cite{\LV, 3.3}, (c) was
proved (in the Weyl group case) by an argument (based on geometry via \cite{\LV, 3.4}) which is not available in 
our case and which we have replaced by the analysis in \S5.)

\head 7. An inversion formula\endhead
\subhead 7.1\endsubhead
In this section we assume that $W$ is finite. 
Let $\huM=\Hom_{\uca}(\uM,\uca)$. For any $w\in\II_*$ we define $\ha'_w\in\huM$ by
$\ha'_w(a'_y)=\d_{y,w}$ for any $y\in\II_*$. Then $\{\ha'_w;w\in\II_*\}$ is an $\uca$-basis of $\huM$. 
We define an $\ufH$-module structure on $\huM$ by $(hf)(m)=f(h^\flat m)$ (with $f\in\huM$, $m\in\uM$, $h\in\ufH$) 
where $h\m h^\flat$ is the algebra antiautomorphism of $\ufH$ such that $T'_s\m T'_s$ for all $s\in S$. (Recall 
that $T'_s=u\i T_s$.) We define a bar operator $\,\bar{}:\huM@>>>\huM$ by $\bar f(m)=\ov{f(\bar m)}$ (with $f\in\huM$, $m\in\uM$); in 
$\ov{f(\bar m)}$ the lower bar is that of $\uM$ and the upper bar is that of $\uca$. We have
$\ov{hf}=\bar h\bar f$ for $f\in\huM$, $h\in\ufH$.

Let $\di:W@>>>W$ be the involution $x\m w_Sx^*w_S=(w_Sxw_S)^*$ which leaves $S$ 
stable. We have $\II_\di=w_S\II_*=\II_*w_S$.
We define the $\uca$-module $\uM_\di$ and its basis $\{b'_z;z\in\II_\di\}$ in terms of
$\di$ in the same way as $\uM$ and its basis $\{a'_w;w\in\II_*\}$ were defined in terms of $*$.
Note that $\uM_\di$ has an $\ufH$-module structure and a bar operator $\,\bar{}:\uM_\di@>>>\uM_\di$
analogous to those of $\uM$.

We define an isomorphism of $\uca$-modules $\Ph:\huM@>>>\uM_\di$ by $\Ph(\ha'_w)=\k(w)b'_{ww_S}$. Here
$\k(w)$ is as in 4.5(a).
Let $h\m h^\dag$ be the algebra automorphism of $\ufH$ such that $T'_s\m-T'_s{}\i$ for any $s\in S$.
We have the following result.

\proclaim{Lemma 7.2} For any $f\in\huM$, $h\in\ufH$ we have $\Ph(hf)=h^\dag\Ph(f)$.
\endproclaim
It is enough to show this when $h$ runs through a set of algebra generators of $\ufH$ and $f$ runs through a
basis of $\huM$. Thus it is enough to show for any $w\in\II_*,s\in S$ that $\Ph(T_s\ha'_w)=-T_s\i\Ph(\ha'_w)$ or
that

(a) $\Ph(T_s\ha'_w)=-\k(w)T_s\i b'_{ww_S}$.
\nl
We write the formulas in 4.1 with $*$ replaced by $\di$ and $a'_w$ replaced by $b'_{ww_S}$:
$$\align&T'_sb'_{ww_S}=b'_{ww_S}+(v+v\i)b'_{sww_S}\text{ if }sw=ws^*<w,\\&
T'_sb'_{ww_S}=(u-1-u\i)b'_{ww_S}+(v-v\i)b'_{sww_S}\text{ if }sw=ws^*>w,\\&
T'_sb'_{ww_S}=b'_{sws^*w_S}\text{ if }sw\ne ws^*<w,\\&
T'_sb'_{ww_S}=(u-u\i)b'_{ww_S}+b'_{sws^*w_S}\text{ if }sw\ne ws^*>w.\endalign$$
Since $T'_s{}\i=T'_s+u\i-u$ we see that
$$\align&-T'_s{}\i b'_{ww_S}=-(u\i+1-u)b'_{ww_S}-(v+v\i)b'_{sww_S}\text{ if $sw=ws^*<w$}\\& 
-T'_s{}\i b'_{ww_S}=b'_{ww_S}-(v-v\i)b'_{sww_S}\text{ if $sw=ws^*>w$}\\& 
-T'_s{}\i b'_{ww_S}=-(u\i-u)b'_{ww_S}-b'_{sws^*w_S}\text{ if $sw\ne ws^*<w$}\\& 
-T'_s{}\i b'_{ww_S}=-b'_{sws^*w_S}\text{ if $sw\ne ws^*>w$}\tag b\endalign$$
Using again the formulas in 4.1 for $T'_sa'_y$ we see that for $y,w\in\II_*$ we have
$$\align&(T'_s\ha'_w)(a_y)=\ha'_w(T'_sa_y)\\&=
\d_{sy=ys^*>y}\d_{y,w}+\d_{sy=ys^*>y}\d_{sy,w}(v+v\i) +\d_{sy=ys^*<y}\d_{y,w}(u-1-u\i)\\&+
\d_{sy=ys^*<y}\d_{sy,w}(v-v\i)+\d_{sy\ne ys^*>y}\d_{sys^*,w}
+\d_{sy\ne ys^*<y}\d_{y,w}(u-u\i)\\&+\d_{sy\ne ys^*<y}\d_{sys^*,w}\\&
=\d_{sw=ws^*>w}\d_{y,w}+\d_{sw=ws^*<w}\d_{y,sw}(v+v\i)+\d_{sw=ws^*<w}\d_{y,w}(u-1-u\i)\\&+
\d_{sw=ws^*>w}\d_{y,sw}(v-v\i)+\d_{sw\ne ws^*<w}\d_{y,sws^*}\\&
+\d_{sw\ne ws^*<w}\d_{y,w}(u-u\i)+\d_{sw\ne ws^*>w}\d_{y,sws^*}\\&
=(\d_{sw=ws^*>w}\ha'_w+\d_{sw=ws^*<w}(v+v\i)\ha'_{sw}+\d_{sw=ws^*<w}(u-1-u\i)\ha'_w\\&+
\d_{sw=ws^*>w}(v-v\i)\ha'_{sw}+\d_{sw\ne ws^*<w}\ha'_{sws^*}\\&
+\d_{sw\ne ws^*<w}(u-u\i)\ha'_w+\d_{sw\ne ws^*>w}\ha'_{sws^*})(a_y).\endalign$$
Since this holds for any $y\in\II_*$ we see that 
$$\align&T'_s\ha'_w=\d_{sw=ws^*>w}\ha'_w+\d_{sw=ws^*<w}(v+v\i)\ha'_{sw}+\d_{sw=ws^*<w}(u-1-u\i)\ha'_w\\&+
\d_{sw=ws^*>w}(v-v\i)\ha'_{sw}+\d_{sw\ne ws^*<w}\ha'_{sws^*}\\&
+\d_{sw\ne ws^*<w}(u-u\i)\ha'_w+\d_{sw\ne ws^*>w}\ha'_{sws^*}.\endalign$$
Thus we have
$$\align&T'_s\ha'_w=\ha'_w+(v-v\i)\ha'_{sw}\text{ if }sw=ws^*>w,\\&
T'_s\ha'_w=(u-1-u\i)\ha'_w+(v+v\i)\ha'_{sw}\text{ if }sw=ws^*<w,\\&
T'_s\ha'_w= \ha'_{sws^*}\text{ if }sw\ne ws^*>w,\\&
T'_s\ha'_w=(u-u\i)\ha'_w+\ha'_{sws^*}\text{ if }sw\ne ws^*<w.\endalign$$
so that
$$\align&\Ph(T'_s\ha'_w)=\k(w)b'_{ww_S}+(v-v\i)\k(sw)b'_{sww_S}\text{ if $sw=ws^*>w$}\\&
\Ph(T'_s\ha'_w)=(u-1-u\i)\k(w)b'_{ww_S}+(v+v\i)\k(sw)b'_{sww_S}\text{ if $sw=ws^*<w$}\\&
\Ph(T'_s\ha'_w)=     \k(sws^*)b'_{sws^*w_S}\text{ if $sw\ne ws^*>w$}\\&
\Ph(T'_s\ha'_w)=(u-u\i)\k(w)b'_w+\k(sws^*)b'_{sws^*w_S}\text{ if $sw\ne ws^*<w$}.\tag c\endalign$$
From (b),(c) we see that to prove (a) we must show:
$$\align&\k(w)b'_{ww_S}+(v-v\i)\k(sw)b'_{sww_S}\\&
=\k(w)b'_{ww_S}-\k(w)(v-v\i)b'_{sww_S}\text{ if $sw=ws^*>w$},\endalign$$
$$\align&(u-1-u\i)\k(w)b'_{ww_S}+(v+v\i)\k(sw)b'_{sww_S}\\&=-\k(w)(u\i+1-u)b'_{ww_S}-\k(w)(v+v\i)b'_{sww_S}     
\text{ if $sw=ws^*<w$},\endalign$$
$$\k(sws^*)b'_{sws^*w_S}=-\k(w)b'_{sws^*w_S}\text{ if $sw\ne ws^*>w$},$$      
$$\align&(u-u\i)\k(w)b'_w+\k(sws^*)b'_{sws^*w_S}\\&
=-\k(w)(u\i-u)b'_{ww_S}-\k(w)b'_{sws^*w_S}\text{ if $sw\ne ws^*<w$}.\endalign$$
This is obvious. The lemma is proved.

\proclaim{Lemma 7.3} We define a map $B:\huM@>>>\huM$ by $B(f)=\Ph\i(\ov{\Ph(f)})$ where the bar refers to 
$\uM_\di$. We have $B(f)=\bar f$ for all $f\in\huM$.
\endproclaim
We show that 

(a) $B(hf)=\bar hB(f)$
\nl
for all $h\in\ufH,f\in\huM$.
This is equivalent to $\Ph\i(\ov{\Ph(hf)})=\bar h\Ph\i(\ov{\Ph(f)})$ or (using 7.2) to
$\ov{h^\dag\Ph(f)}=\Ph(\bar h\Ph\i(\ov{\Ph(f)}))$ or (using 7.2) to
$\ov{h^\dag}(\ov{\Ph(f)})=(\bar h)^\dag \Ph(\Ph\i(\ov{\Ph(f)}))$; it remains to use that
$\ov{h^\dag}=(\bar h)^\dag$.

Next we show that

(b) $B(\ha'_{w_S})=\ha'_{w_S}$.
\nl
Indeed the left hand side is 
$$\Ph\i(\ov{\Ph(\ha'_{w_S})})=\Ph\i(\ov{\k(w_S)b'_1})=\k(w_S)\Ph\i(b'_1)=\ha'_{w_S}$$
as required. (We have used that $\ov{b'_1}=b'_1$ in $\uM_\di$.) Next we show:

(c) $\ov{\ha'_{w_S}}=\ha'_{w_S}$.
\nl
Indeed for $y\in\II_*$ we have 
$$\ov{\ha'_{w_S}}(a'_y)=\ov{\ha'_{w_S}(\ov{a'_y})}
=\ov{\ha'_{w_S}(\sum_{x\in\II_*;x\le y}\bar r_{x,y}a'_x)}=\ov{\bar r_{w_S,w_S}\d_{y,w_S}}=\d_{y,w_S}
=\ha'_{w_S}(a'_y)$$
(we use that $r_{w_S,w_S}=1$). This proves (c).

Since $\ov{hf}=\bar h\bar f$ for all $h\in\ufH,f\in\huM$ we see (using (a),(b),(c)) that
the map $f\m\ov{B(f)}$ from $\huM$ into itself is $\ufH$-linear and carries $\ha'_{w_S}$ to itself.
This implies that this map is the identity. (It is enough to show that $\ha'_{w_S}$ generates the $\ufH$-module
$\huM$ after extending scalars to $\QQ(v)$. Using 7.2 it is enough to show that $b'_1$ generates the
$\ufH$-module $\uM_\di$ after extending scalars to $\QQ(v)$. This is known from 2.11.)
We see that $f=\ov{B(f)}$ for all $f\in\huM$. Applying $\,\bar{}$ to both sides (an involution of $\huM$) we
deduce that $\bar f=B(f)$ for all $f\in\huM$. The lemma is proved.

\subhead 7.4\endsubhead
Recall that $\ov{a'_w}=\sum_{y\in\II_*;y\le w}\ov{r_{y,w}}a'_y$ for $w\in\II_*$. 
The analogous equality in $\uM_\di$ is
$$\ov{b'_z}=\sum_{x\in\II_\di;x\le z}\ov{r^\di_{x,z}}b'_x\text{ for }x\in\II_\di.\tag a$$
Here $r^\di_{x,z}\in\uca$. We have the following result.

\proclaim{Proposition 7.5} Let $y,w\in\II_*$ be such that $y\le w$. We have 
$$\ov{r_{y,w}}=\k(y)\k(w)r^\di_{ww_S,yw_S}.$$
\endproclaim
We show that for any $y\in\II_*$ we have 
$$\ov{\ha'_y}=\sum_{w\in\II_*;y\le w}r_{y,w}\ha'_w.\tag a$$
Indeed for any $x\in\II_*$ we have 
$$\align&\ov{\ha'_y}(a'_x)=\ov{\ha'_y(\ov{a'_x})}=\ov{\ha'_y(\sum_{x'\in\II_*;x'\le x}\bar r_{x',x}a'_{x'})}
=\ov{\d_{y\le x}\bar r_{y,x}}=\d_{y\le x}r_{y,x}\\&=\sum_{w\in\II_*;y\le w}r_{y,w}\ha'_w(a'_x).\endalign$$
Using (a) and 7.3 we see that for any $y\in\II_*$ we have 
$$\Ph\i(\ov{\Ph(\ha'_y)})=\sum_{w\in\II_*;y\le w}r_{y,w}\ha'_w.$$
It follows that $\ov{\Ph(\ha'_y)}=\sum_{w\in\II_*;y\le w}r_{y,w}\Ph(\ha'_w)$ that is,
$$\ov{\k(y)b'_{yw_S}}=\sum_{w\in\II_*;y\le w}r_{y,w}\k(w)b'_{ww_S}.$$
Using 7.4(a) to compute the left hand side we obtain
$$\k(y)\sum_{w\in\II_*;ww_S\le yw_S}\ov{r^\di_{ww_S,yw_S}}b'_{ww_S}
=\sum_{w\in\II_*;y\le w}r_{y,w}\k(w)b'_{ww_S}.$$
Hence for any $w\in\II_*$ such that $y\le w$ we have $r_{y,w}\k(w)=\k(y)\ov{r^\di_{ww_S,yw_S}}$. The
proposition follows.

\subhead 7.6\endsubhead
Recall that for $y,w\in\II_*$, $y\le w$ we have $P^\pm_{y,w}=v^{l(w)-l(y)}\p_{y,w}$ where 
$\p_{y,w}\in\uca$ satisfies $\p_{w,w}=1$, $\p_{y,w}\in\uca_{<0}$ if $y<w$ and
$$\ov{\p_{y,w}}=\sum_{t\in\II_*;y\le t\le w}r_{y,t}\p_{t,w}.\tag a$$
Replacing $*$ by $\di$ in the definition of $P^\pm_{y,w}$ we obtain polynomials
$P^{\s,\di}_{x,z}\in\ZZ[u]$ ($x,z\in\II_\di,x\le z$) such that 
$P^{\s,\di}_{x,z}=v^{l(z)-l(x)}\p^\di_{x,z}$ where 
$\p^\di_{x,z}\in\uca$ satisfies $\p^\di_{z,z}=1$, $\p^\di_{x,z}\in\uca_{<0}$ if $x<z$ and
$$\ov{\p^\di_{x,z}}=\sum_{t'\in\II_\di;x\le t'\le z}r^\di_{x,t'}\p^\di_{t',z}.\tag b$$
The following inversion formula (and its proof) is in the same spirit as \cite{\KL, 3.1} (see also \cite{\VO}).

\proclaim{Theorem 7.7}For any $y,w\in\II_*$ such that $y\le w$ we have 
$$\sum_{t\in\II_*;y\le t\le w}\k(y)\k(t)P^\pm_{y,t}P^{\pm,\di}_{ww_S,tw_S}=\d_{y,w}.$$
\endproclaim
The last equality is equivalent to
$$\sum_{t\in\II_*;y\le t\le w}\k(y)\k(t)\p_{y,t}\p^\di_{ww_S,tw_S}=\d_{y,w}.\tag a$$
Let $M_{y,w}$ be the left hand side of (a). When $y=w$ we have $M_{y,w}=1$. Thus, we may assume that $y<w$ and 
that $M_{y',w'}=0$ for all $y',w'\in\II_*$ such that $y'<w'$, $l(w')-l(y')<l(w)-l(y)$. Using 7.6(a),(b) we have
$$\align&M_{y,w}=\sum_{t\in\II_*;y\le t\le w}\k(y)\k(t)\sum_{x,x'\in\II_*;y\le x\le t\le x'\le w}
\ov{r_{y,x}}\ov{p_{x,t}}\ov{r^\di_{ww_S,x'w_S}}\ov{p^\di_{x'w_S,tw_S}}\\&
=\sum_{x,x'\in\II_*;y\le x\le x'\le w}\k(y)\k(x)\ov{r_{y,x}}\ov{r^\di_{ww_S,x'w_S}}\ov{M_{x,x'}}.\endalign$$
The only $x,x'$ which can contribute to the last sum satisfy $x=x'$ or $x=y,x'=w$. Thus
$$M_{y,w}=\sum_{x\in\II_*;y\le x\le w}\k(y)\k(x)\ov{r_{y,x}}\ov{r^\di_{ww_S,xw_S}}+\ov{M_{y,w}}.$$
(We have used 4.8(a).) Using 7.5 we see that the last sum over $x$ is equal to
$$\k(y)\k(w)\sum_{x\in\II_*;y\le x\le w}\ov{r_{y,x}}r_{x,w}=0,$$
see 4.6(a). Thus we have $M_{y,w}=\ov{M_{y,w}}$. Since $M_{y,w}\in\uca_{<0}$, this forces $M_{y,w}=0$. The 
theorem is proved.

\head 8. A $(-u)$ analogue of weight multiplicities?\endhead
\subhead 8.1\endsubhead
In this section we assume that $W$ is an irreducible affine Weyl group.
An element $x\in W$ is said to be a translation if its $W$-conjugacy class is finite. The set of translations is
a normal subgroup $\ct$ of $W$ of finite index. We fix an element $s_0\in S$ such that, setting $K=S-\{s_0\}$,
the obvious map $W_K@>>>W/\ct$ is an isomorphism. (Such an $s_0$ exists.) We assume that $*$ is the automorphism 
of $W$ such that $x\m w_Kxw_K$ for all $x\in W_K$ and $y\m w_Ky\i w_K$ for any $y\in\ct$ (this automorphism maps 
$s_0$ to $s_0$ hence it maps $S$ onto itself). We have $K^*=K$.

\proclaim{Proposition 8.2} If $x$ is an element of $W$ which has maximal length in its $(W_K,W_K)$ double coset
$\Om$ then $x^*=x\i$.
\endproclaim
Note that $\ct_\Om:=\Om\cap\ct$ is a single $W$-conjugacy class. If $y\in\ct_\Om$ then 
$y^{*-1}=w_Kyw_K\in\ct_\Om$. Thus $w\m w^{*-1}$ maps some element of $\Om$ to an element of $\Om$. Hence it maps 
$\Om$ onto itself. Since it is length preserving it maps $x$ to itself.

\subhead 8.3\endsubhead
Let $\Om,\Om'$ be two $(W_K,W_K)$-double cosets in $W$ such that $\Om'\le\Om$. As in 5.1, let $d_\Om$ (resp. 
$d_{\Om'}$) be the longest element in $\Om$ (resp. $\Om'$). Let $P_{d_{\Om'},d_\Om}\in\ZZ[u]$ be the 
polynomial attached in \cite{\KL} to the elements $d_{\Om'},d_\Om$ of the Coxeter group $W$. Let $G$ be a 
simple adjoint group over $\CC$ for which $W$ is the associated affine Weyl group so that $\ct$ is the lattice of
weights of a maximal torus of $G$. Let $V_\Om$ be the (finite dimensional) irreducible rational representation of
$G$ whose extremal weights form the set $\ct_\Om$. Let $N_{\Om',\Om}$ be the multiplicity of a weight in 
$\ct_{\Om'}$ in the representation $V_\Om$. Now $P_{d_{\Om'},d_\Om}$ is the $u$-analogue (in the sense of 
\cite{\LSI}) of the weight multiplicity $N_{\Om',\Om}$; in particular, according to \cite{\LSI}, we have
$$N_{\Om',\Om}=P_{d_{\Om'},d_\Om}|_{u=1}.$$
We have the following 
\proclaim{Conjecture 8.4} $P^\pm_{d_{\Om'},d_\Om}(u)=P_{d_{\Om'},d_\Om}(-u)$.
\endproclaim

\subhead 8.5\endsubhead
Now assume that $\Om$ (resp. $\Om'$) is the $(W_K,W_K)$-double coset that contains $s_0$ (resp. the unit 
element). Let $e_1\le e_2\le\do\le e_n$ be the exponents of $W_K$ (recall that $e_1=1$). The following result 
supports the conjecture in 8.4.

\proclaim{Proposition 8.6} In the setup of 8.5, assume that $W_K$ is simply laced. We have:

(a) $A_{d_\Om}=v^{-l(d_\Om)}a_\Om+(-1)^{e_n}\sum_{j\in[1,n]}(-u)^{-e_j}v^{-l(d_{\Om'})}a_{\Om'};$

(b) $P_{d_{\Om'},d_\Om}(u)=\sum_{j\in[1,n]}u^{e_j-1}$; 

(c) $P^\pm_{d_{\Om'},d_\Om}(u)=\sum_{j\in[1,n]}(-u)^{e_j-1}.$
\endproclaim
We prove (a). It is enough to show that 
$$v^{-l(d_\Om)}a_\Om+(-1)^{e_n}\sum_{j\in[1,n]}(-u)^{-e_j}v^{-l(d_{\Om'})}a_{\Om'}$$
is fixed by $\,\bar{}$. Let $H=K\cap s_0Ks_0$. We have $H=H^*$ and $W_H$ is contained in the centralizer of $s_0$.
Let $\t:W_H@>>>W_H$ be the automorphism $y\m s_0y^*s_0=y^*$. We have $d_{\Om'}=w_K$, $d_\Om=w_Kw_Hs_0w_K$,
$l(d_\Om)=2l(w_K)-l(w_H)+1$ and we must show that
$$v^{-2l(w_K)+l(w_H)-1}a_\Om+(-1)^{e_n}\sum_{j\in[1,n]}(-u)^{-e_j}v^{-l(w_K)}a_{\Om'}\text{ is fixed by }\,\bar{}.
\tag d$$
Let $\SS=\sum_{x\in W_K}T_x\in\ufH$. Using 5.10(a) we see that 
$$\SS(a_{s_0}+a_1)=\PP_{H,*}a_\Om+\PP_{K,*}a_{\Om'}.$$
Hence
$$\align&v^{-2l(w_K)}\SS(v\i(a_{s_0}+a_\em))\\&=v^{-l(w_H)}\PP_{H,*}v^{-2l(w_K)+l(w_H)-1}a_\Om
+v^{-l(w_K)-1}\PP_{K,*} v^{-l(w_K)}a_{\Om'}.\endalign$$
Since $v^{-2l(w_K)}\SS$ and $v\i(a_{s_0}+a_1)$ are fixed by $\,\bar{}$, we see that
that the left hand side of the last equality is fixed by $\,\bar{}$, hence
$$v^{-l(w_H)}\PP_{H,*}v^{-2l(w_K)+l(w_H)-1}a_\Om+v^{-l(w_K)-1}\PP_{K,*}v^{-l(w_K)}a_{\Om'}$$
is fixed by $\,\bar{}$. Since $v^{-l(w_H)}\PP_{H,*}$ is fixed by $\,\bar{}$ and divides $\PP_{K,*}$, we see that 
$$v^{-2l(w_K)+l(w_H)-1}a_\Om+v^{-l(w_K)+l(w_H)-1}\PP_{K,*}\PP_{H,*}\i v^{-l(w_K)}a_{\Om'}$$
is fixed by $\,\bar{}$. Hence to prove (d) it is enough to show that
$$v^{-l(w_K)+l(w_H)-1}\PP_{K,*}\PP_{H,*}\i v^{-l(w_K)}a_{\Om'}-
(-1)^{e_n}\sum_{j\in[1,n]}(-u)^{-e_j}v^{-l(w_K)}a_{\Om'}$$
is fixed by $\,\bar{}$. Now $v^{-l(w_K)}a_{\Om'}$ is fixed by $\,\bar{}$, see 5.11(a). Hence it is enough to show 
that
$$v^{-l(w_K)+l(w_H)-1}\PP_{K,*}\PP_{H,*}\i-(-1)^{e_n}\sum_{j\in[1,n]}(-u)^{-e_j}\text{ is fixed by }\,\bar{}.$$
This is verified by direct computation in each case. This completes the proof of (a). Now (c) follows from (a) 
using the equality $l(w_Kw_Hs_0w_K)-l(w_K)=2e_n$ and the known symmetry property of exponents; (b) follows from 
\cite{\LSI}.

\subhead 8.7\endsubhead
In this subsection we assume that $W_K$ is of type $A_2$ with $K=\{s_1,s_2\}$. Note that $s_1^*=s_2$, 
$s_2^*=s_1$. We write $i_1i_2\do$ instead of $s_{i_1}s_{i_2}\do$ (the indices are in $\{0,1,2\}$). Let 
$\Om_1,\Om_2,\Om_3,\Om_4,\Om_5$ be the $(W_K,W_K)$ double coset of $01210,0120,0210,0$ and unit element
respectively. We have $d_{\Om_1}=1210120121$, $d_{\Om_2}=121012012$, $d_{\Om_3}=121021021$, 
$d_{\Om_4}=1210121$,  $d_{\Om_5}=121$. A direct computation shows that
$$A_{d_{\Om_1}}=v^{-11}(a_{\Om_1}+a_{\Om_2}+a_{\Om_3}+(1-u)a_{\Om_4}+(1-u+u^2)a_{\Om_5}).$$
This provides further evidence for the conjecture in 8.4.

\subhead 8.8\endsubhead
In this subsection we assume that $K=\{s_1,s_2\}$ with $s_1s_2$ of order $4$ and with $s_0s_2=s_2s_0$, $s_0s_1$ 
of order $4$. Note that $x^*=x$ for all $x\in W$. Let $\Om_1,\Om_2,\Om_3$ be the $(W_K,W_K)$ double coset of 
$s_0s_1s_0,s_0$ and unit element respectively. We have $d_{\Om_1}=1212010212$, $d_{\Om_2}=12120121$, 
$d_{\Om_3}=1212$ (notation as in 8.7). A direct computation shows that
$$A_{d_{\Om_1}}=v^{-10}(a_{\Om_1}+a_{\Om_2}+(1+u^2)a_{\Om_3}).$$
This provides further evidence for the conjecture in 8.4.

\head 9. Reduction modulo $2$\endhead
\subhead 9.1\endsubhead
Let $\ca_2=\ca/2\ca=(\ZZ/2)[u,u\i]$, $\uca_2=\uca/2\uca=(\ZZ/2)[v,v\i]$. We regard $\ca_2$ as a subring of 
$\uca_2$ by setting $u=v^2$. Let $\fH_2=\fH/2\fH$; this is naturally an $\ca_2$-algebra with $\ca_2$-basis 
$(T_x)_{x\in W}$ inherited from $\fH$ and with a bar operator $\bar{}:\fH_2@>>>\fH_2$ inherited from that
$\fH$. Let $M_2=\ca_2\ot_\ca M=M/2M$. This has a $\fH_2$-module structure and a bar operator
$\bar{}:M_2@>>>M_2$ inherited from $M$. It has an $\ca_2$-basis $(a_w)_{w\in\II_*}$ inherited from $M$.
In this section we give an alternative construction of the $\fH_2$-module structure on $M_2$ and its bar operator.

Let $\ch$ be the free $\uca$-module with basis 
$(t_w)_{w\in W}$ with the unique $\uca$-algebra structure with unit $t_1$ such that 

$t_wt_{w'}=t_{ww'}$ if $l(ww')=l(w)+l(w')$ and

$(t_s+1)(t_s-v^2)=0$ for all $s\in S$.
\nl
Let $\,\bar{}\,:\ch@>>>\ch$ be the unique ring involution such that $\ov{v^nt_x}=v^{-n}t_{x\i}\i$ for any
$x\in W,n\in\ZZ$ (see \cite{\KL}). Let $\ch_2=\ch/2\ch$; this is naturally an $\uca_2$-algebra with 
$\uca_2$-basis $(t_x)_{x\in W}$ inherited from $\ch$ and with a bar operator $\bar{}:\ch_2@>>>\ch_2$ inherited 
from that of $\ch$. Let $h\m h^\sp$ be the unique algebra antiautomorphism of $\ch$ such that 
$t_w\m t_{w^{*-1}}$. (It is an involution.)

We have $\ch_2=\ch'_2\op\ch''_2$ where $\ch'_2$ (resp. $\ch''_2$) is the $\uca$-submodule of $\ch_2$ spanned by
$\{t_w;w\in\II_*\}$ (resp. $\{t_w;w\in W-\II_*\}$). Let $\p:\ch_2@>>>\ch'_2$ be the projection on the first summand.
Note that for $\x'\in\ch_2$ we have

(a) $\x'{}^\sp=\x'$ if and only if $\x'=\x'_1+\x'_2+\x'_2{}^\sp$ where $\x'_1\in\ch'_2,\x'_2\in\ch_2$.

(b) $\p(\x'{}^\sp)=\p(\x')$.

\proclaim{Lemma 9.2} The map $\ch_2\T\ch'_2@>>>\ch'_2$, $(h,\x)\m h\circ\x=\p(h\x h^\sp)$
defines an $\ch_2$-module structure on the abelian group $\ch'_2$.
\endproclaim
Let $h,h'\in\ch_2,\x\in\ch'_2$. We first show that $(h+h')\circ\x=h\circ\x+h'\circ\x$ or that 
$\p((h+h')\x(h+h')^\sp)=\p(h\x h^\sp)+\p(h'\x h'{}^\sp)$. It is enough to show that
$\p(h\x h'{}^\sp)=\p(h'\x h^\sp)$. This follows from 9.1(b) since
$(h'\x h^\sp)^\sp=h\x^\sp h'{}^\sp=h\x h'{}^\sp$.

We next show that $(hh')\circ\x=h\circ(h'\circ\x)$ or that 
$\p(hh'\x h'{}^\sp h^\sp)= \p(h\p(h'\x h'{}^\sp)h^\sp)$. Setting $\x'=h'\x h'{}^\sp$ we see that we
must show that $\p(h\x'h^\sp)= \p(h\p(\x')h^\sp)$. 
Setting $\et=\x'-\p(\x')$ we are reduced to showing that $\p(h\et h^\sp)=0$.
Since $\x\in\ch'_2$ we have $\x^\sp=\x$. Hence $\x'{}^\sp=(h'{}^\sp)^\sp\x^\sp h'{}^\sp=h'\x h'{}^\sp$
so that $\x'{}^\sp=\x'$. We write $\x'=\x'_1+\x'_2+\x'_2{}^\sp$ as in 9.1(a). Then $\p(\x')=\x'_1$ and
$\et=\x'_2+\x'_2{}^\sp$. We have $h\et h^\sp=h\x'_2h^\sp+h\x'_2{}^\sp h^\sp=\z+\z^\sp$ where $\z=h\x'_2h^\sp$.
Thus $\p(h\et h^\sp)=\p(\z+\z^\sp)=0$ (see 9.1(b)). Clearly we have $1\circ\x=\x$. The lemma is proved.

\subhead 9.3\endsubhead
Consider the group isomorphism $\ps:\ch_2@>\sim>>\fH_2$ such that $v^nt_w\m u^nT_w$ for any $n\in\ZZ,w\in W$. 
This is a ring isomorphism satisfying $\ps(fh)=f^2\ps(h)$ for all $f\in\uca_2,h\in\ch_2$ (we have 
$f^2\in\ca_2$). Using now 9.2 we see that:

(a) {\it The map $\fH_2\T\ch'_2@>>>\ch'_2$, $(h,\x)\m h\odot\x:=\p(\ps\i(h)\x(\ps\i(h))^\sp)$ defines an 
$\fH_2$-module structure on the abelian group $\ch'_2$.}
\nl
Note that the $\fH_2$-module structure on $\ch'_2$ given in (a) is compatible with the $\ca$-module structure on
$\ch'_2$. Indeed if $f\in\ca_2$ and $f'\in\uca_2$ is such that $f'{}^2=f$ then $f$ acts in the $\fH_2$-module
structure in (a) by $\x\m f'\x f'=f'{}^2\x=f\x$.

\subhead 9.4\endsubhead
Let $s\in S,w\in\II_*$. The equation in this subsection take place in $\ch_2$. If $sw=ws^*>w$ we have
$$T_s\odot t_w=\p(t_s t_wt_{s^*})=\p(t_{sw}t_{s^*})=\p((u-1)t_{sw}+ut_w)=ut_w+(u+1)t_{sw}.$$
If $sw=ws^*<w$ we have
$$\align&T_s\odot t_w=\p(t_s t_wt_{s^*})=\p(((u-1)t_w+ut_{sw})t_{s^*})\\&
=\p((u-1)^2t_w+(u-1)ut_{ws^*}+ut_w)=(u^2-u-1)t_w+(u^2-u)t_{sw}.\endalign$$
If $sw\ne ws^*>w$ we have
$$T_s\odot t_w=\p(t_s t_wt_{s^*})=\p(t_{sws^*})=t_{sws^*}.$$
If $sw\ne ws^*<w$ we have
$$\align&T_s\odot t_w=\p(t_s t_wt_{s^*})=\p(((u-1)t_w+ut_{sw})t_{s^*})\\&
=\p((u-1)^2t_w+(u-1)ut_{ws^*}+u(u-1)t_{sw}+u^2t_{sws^*})=(u^2-1)t_w+u^2t_{sws^*}.\endalign$$
(We have used that $\p(t_{ws^*})=\p(t_{sw})$ which follows from 9.1(b).)
From these formulas we see that 

(a) {\it the isomorphism of $\ca_2$-modules $\ch'_2@>\sim>>M_2$ given by $t_w\m a_w$ ($w\in\II_*$) is compatible
with the $\fH_2$-module structures.}

\subhead 9.5\endsubhead
For $w\in W$ we set $\ov{t_w}=\sum_{y\in W;y\le w}\ov{\r_{y,w}}v^{-l(w)-l(y)}t_y$ where $\r_{y,w}\in\uca$ 
satisfies $\r_{w,w}=1$. For $y\in W,y\not\le w$ we set $\r_{y,w}=0$.

For $x,y\in W,s\in S$ such that $sy>y$ we have

(i) $\r_{x,sy}=\r_{sx,y}$ if $sx<x$,

(ii) $\r_{x,sy}=\r_{sx,y}+(v-v\i)\r_{x,y}$ if $sx>x$.

For $x,y\in W,s\in S$ such that $ys>y$ we have

(iii) $\r_{x,ys}=\r_{xs,y}$ if $xs<x$,

(iv) $\r_{x,ys}=\r_{xs,y}+(v-v\i)\r_{x,y}$ if $xs>x$.

Note that (iii),(iv) follow from (i),(ii) using

(v) $\r_{z,w}=\r_{z^{*-1},w^{*-1}}$ for any $z,w\in W$.

\subhead 9.6\endsubhead
If $f,f'\in\uca$ we write $f\equiv f'$ if $f,f'$ have the same image under the obvious ring homomorphism 
$\uca@>>>\uca_2$. We have the following result.

\proclaim{Proposition 9.7} For any $y,w\in\II_*$ we have $r_{y,w}\equiv\r_{y,w}$.
\endproclaim
Since the formulas 4.2(a),(b) together with $r_{x,1}=\d_{x,1}$ define uniquely $r_{x,y}$ for any $x,y\in\II_*$ 
and since $\r_{x,1}=\d_{x,1}$ for any $x$, it is enough to show that the equations 4.2(a),(b) remain valid if each
$r$ is replaced by $\r$ and each $=$ is replaced by $\equiv$. 

Assume first that $sy=ys^*>y$ and $x\in\II_*$.

If $sx=xs^*>x$ we have
$$\align&(v+v\i)\r_{x,sy}-(\r_{sx,y}(v\i-v)-(u-u\i)\r_{x,y})\\&
\equiv(v+v\i)(\r_{x,sy}-\r_{sx,y}-(v-v\i)\r_{x,y})=0.\endalign$$
(The $=$ follows from 9.5(ii).)

If $sx=xs^*<x$ we have
$$(v+v\i)\r_{x,sy}-(-2\r_{x,y}+\r_{sx,y}(v+v\i))\equiv(v+v\i)(\r_{x,sy}-\r_{sx,y})=0.$$
(The $=$ follows from 9.5(i).)

If $sx\ne xs^*>x$ we have
$$\align&(v+v\i)\r_{x,sy}-(\r_{sxs^*,y}+(u-1-u\i)\r_{x,y})\\&=
(v+v\i)\r_{sx,y}+(u-u\i)\r_{x,y}-\r_{sxs^*,y}-(u-1-u\i)\r_{x,y}\equiv\\&
(v-v\i)\r_{sx,y}-\r_{x,y}+\r_{sxs^*,y}=\r_{sx,ys^*}-\r_{x,y}=0.\endalign$$
(The first, second and third $=$ follow from 9.5(ii),(iv),(iii).)

If $sx\ne xs^*<x$ we have
$$\align&(v+v\i)\r_{x,sy}-(-\r_{x,y}+\r_{sxs^*,y})=(v+v\i)\r_{sx,y}-(-\r_{x,y}+\r_{sxs^*,y})\equiv\\&
(v-v\i)\r_{sx,y}+\r_{x,y}-\r_{sxs^*,y}=\r_{sx,sy}-\r_{sxs^*,y}=\r_{sx,sy}-\r_{sx,ys^*}=0.\endalign$$
(The first, second and third $=$ follow from 9.5(i),(ii),(iii).)

Next we assume that $sy\ne ys^*>y$ and $x\in\II_*$.

If $sx=xs^*>x$ we have
$$\align&\r_{x,sys^*}-(\r_{sx,y}(v\i-v)+(u+1-u\i)\r_{x,y})\\&=\r_{sx,ys^*}+(v-v\i)\r_{x,ys^*}-\r_{sx,y}(v\i-v)
-(u+1-u\i)\r_{x,y}\\&=\r_{x,y}+(v-v\i)\r_{x,ys^*}-\r_{xs^*,y}(v\i-v)-(u+1-u\i)\r_{x,y}\\&=
\r_{x,y}+(v-v\i)\r_{xs^*,y}+(v-v\i)^2\r_{x,y}-\r_{xs^*,y}(v\i-v)\\&-(u+1-u\i)\r_{x,y}\equiv0.\endalign$$
(The first, second and third $=$ follow from 9.5(ii),(iv),(iv).)

If $sx=xs^*<x$ we have
$$\align&\r_{x,sys^*}-(\r_{sx,y}(v+v\i)-\r_{x,y})=\r_{sx,sy}-(\r_{sx,y}(v+v\i)-\r_{x,y})\equiv\\&
\r_{sx,sy}-(\r_{sx,y}(v-v\i)+\r_{x,y})=0,\endalign$$
(The first and second $=$ follow from 9.5(i),(ii.)

If $sx\ne xs^*>x$ we have
$$\align&\r_{x,sys^*}-(\r_{sxs^*,y}+(u-u\i)\r_{x,y})\\&=\r_{xs^*,sy}+(v-v\i)\r_{x,sy}-\r_{sxs^*,y}-(u-u\i)
\r_{x,y}\\&=\r_{sxs^*,y}+(v-v\i)\r_{xs^*,y}+(v-v\i)\r_{sx,y}+(v-v\i)^2\r_{x,y}\\&-\r_{sxs^*,y}-(u-u\i)\r_{x,y}
\equiv(v-v\i)
(\r_{xs^*,y}-\r_{sx,y})\\&=(v-v\i)(\r_{(xs^*)^{*-1},y^{*-1}}-\r_{sx,y})=(v-v\i)(\r_{sx,y}-\r_{sx,y})=0.\endalign$$
(The first, second, and third $=$ follow from 9.5(iv),(ii),(v).)

If $sx\ne xs^*<x$ we have
$$\r_{x,sys^*}-\r_{sxs^*,y}=\r_{xs^*,ys^*}-\r_{sxs^*,y}=0.$$
(The first and second $=$ follow from 9.5(iii),(i).)

Thus the equations 4.2(a),(b) with each $r$ replaced by $\r$ and each $=$ replaced by $\equiv$ are verified. The 
proposition is proved.

\subhead 9.8\endsubhead
We define a group homomorphism $B:\ch'_2@>>>\ch'_2$ by $\x\m\p(\ov{\x})$. From 9.7 we see that 

(a) {\it under the isomorphism 9.4(a) the map $B:\ch'_2@>>>\ch'_2$ corresponds to the map $\,\bar{}:M_2@>>>M_2$.}
\nl
We now give an alternative proof of (a). Using 0.2(b) and 9.4(a) we see that it is enough to show that for any 
$w\in\II_*$ we have $\p(t_{w\i}\i)=T_{w\i}\i\odot t_{w\i}$ in $\ch'_2$. Since $\ps$ in 9.3 is a ring isomorphism,
we have $\ps(t_{w\i}\i)=T_{w\i}\i$ hence 
$$\align&T_{w\i}\i\odot t_{w\i}=\p(\ps\i(T_{w\i}\i)t_{w\i}(\ps\i(T_{w\i}\i))^\sp)\\&=
\p(t_{w\i}\i t_{w\i}(t_{w\i}\i)^\sp)=\p(t_{w\i}\i t_{w\i}t_{w^*}\i)=\p(t_{w\i}\i t_{w\i}t_{w\i}\i)=
\p(t_{w\i}),\endalign$$
as required.

\subhead 9.9\endsubhead
For $y,w\in W$ let $P_{y,w}\in\ZZ[u]$ be the polynomials defined in \cite{\KL, 1.1}. (When 
$y\not\le w$ we set $P_{y,w}=0$.) We set $p_{y,w}=v^{-l(w)+l(y)}P_{y,w}\in\uca$. Note that $p_{w,w}=1$
and $p_{y,w}=0$ if $y\not\le w$. We have $p_{y,w}\in\ca_{<0}$ if $y<w$ and

(i) $\ov{p_{x,w}}=\sum_{y\in W;x\le y\le w}r_{x,y}p_{y,w}$ if $x\le w$,

(ii) $p_{x^{*-1},w^{*-1}}=p_{x,w}$, if $x\le w$.
\nl 
We have the following result which, in the special case where $W$ is a Weyl group or an affine Weyl group, can be
deduced from the last sentence in the first paragraph of \cite{\LV}.

\proclaim{Theorem 9.10} For any $x,w\in\II_*$ such that $x\le w$ we have $P^\pm_{x,w}\equiv P_{x,w}$ (with 
$\equiv$ as in 9.6).
\endproclaim
It is enough to show that $\p_{x,w}\equiv p_{x,w}$. We can assume that $x<w$ and that the result is known when 
$x$ is replaced by $x'\in\II_*$ with $x<x'\le w$. Using 9.9(i) and the definition of $\p_{x,w}$ we have
$$\ov{p_{x,w}}-\ov{\p_{x,w}}
=\sum_{y\in W;x\le y\le w}r_{x,y}p_{y,w}-\sum_{y\in\II_*;x\le y\le w}\r_{x,y}\p_{y,w}.$$
Using 9.7 and the induction hypothesis we see that the last sum is $\equiv$ to
$$\align&p_{x,w}-\p_{x,w}+\sum_{y\in W;x<y\le w}r_{x,y}p_{y,w}-\sum_{y\in\II_*;x<y\le w}r_{x,y}p_{y,w}\\&=
p_{x,w}-\p_{x,w}+\sum_{y\in W;y\ne y^{*-1},x<y\le w}r_{x,y}p_{y,w}.\endalign$$
In the last sum the terms corresponding to $y$ and $y^{*-1}$ cancel out (after reduction $\mod2$) since
$$r_{x,y^{*-1}}p_{y^{*-1},w}=r_{x^{*-1},y}p_{y,w^{*-1}}=r_{x,y}p_{y,w}.$$
(We use 9.5(v), 9.9(ii).) We see that 
$$\ov{p_{x,w}}-\ov{\p_{x,w}}\equiv p_{x,w}-\p_{x,w}.$$
After reduction $\mod2$ the right hand side is in $v\i(\ZZ/2)[v\i]$ and the left hand side is in $v(\ZZ/2)[v]$; 
hence both sides are zero in $(\ZZ/2)[v,v\i]$. This completes the proof.

\subhead 9.11\endsubhead
For $x,w\in\II_*$ such that $x\le w$ we set $P^+_{x,w}=(1/2)(P_{x,w}+P^\pm_{x,w})$, 
$P^-_{x,w}=(1/2)(P_{x,w}-P^\pm_{x,w})$. From 9.10 we see that $P^+_{x,w}\in\ZZ[u]$, $P^-_{x,w}\in\ZZ[u]$.

\proclaim{Conjecture 9.12} We have $P^+_{x,w}\in\NN[u]$, $P^-_{x,w}\in\NN[u]$.
\endproclaim
This is a refinement of the conjecture in \cite{\KL} that $P_{x,w}\in\NN[u]$ for any $x\le w$ in $W$. In the case
where $W$ is a Weyl group or an affine Weyl group, the (refined) conjecture holds by results of \cite{\LV}.

\enddocument